\documentclass[12pt,oneside,reqno]{amsart}
\usepackage{}
\usepackage{amssymb}

\usepackage{bbm}
\usepackage{cases}
\usepackage{amsmath}
\usepackage{graphicx}
\usepackage{mathrsfs}
\usepackage{stmaryrd}
\usepackage{color}
\usepackage{soul}
\usepackage[dvipsnames]{xcolor}
\usepackage{amsfonts}
\usepackage{enumerate,amsmath,amssymb,amsthm}

\numberwithin{equation}{section}

\newcommand{\be}{\begin{eqnarray}}
\newcommand{\mE}{\end{eqnarray}}
\newcommand{\ce}{\begin{eqnarray*}}
\newcommand{\de}{\end{eqnarray*}}
\newtheorem{theorem}{Theorem}[section]
\newtheorem{lemma}[theorem]{Lemma}
\newtheorem{remark}[theorem]{Remark}
\newtheorem{definition}[theorem]{Definition}
\newtheorem{proposition}[theorem]{Proposition}
\newtheorem{example}[theorem]{Example}
\newtheorem{corollary}[theorem]{Corollary}

\def\eps{\varepsilon}

\def\p{\partial}

\def\[{{\Big[}}
\def\]{{\Big]}}
\def\<{{\langle}}
\def\>{{\rangle}}
\def\({{\Big(}}
\def\){{\Big)}}

\def\bx{{\mathbf{x}}}

\def\dif{{\mathord{{\rm d}}}}

\def\no{\nonumber}
\def\={&\!\!=\!\!&}
\def\bt{\begin{theorem}}
\def\et{\end{theorem}}
\def\bl{\begin{lemma}}
\def\el{\end{lemma}}
\def\br{\begin{remark}}
\def\er{\end{remark}}

\def\bd{\begin{definition}}
\def\ed{\end{definition}}
\def\bp{\begin{proposition}}
\def\ep{\end{proposition}}
\def\bc{\begin{corollary}}
\def\ec{\end{corollary}}
\def\bx{\begin{example}}
\def\ex{\end{example}}

\def\cI{{\mathcal I}}
\def\cJ{{\mathcal J}}
\def\cK{{\mathcal K}}

\def\cQ{{\mathcal Q}}

\def\mE{{\mathbb E}}

\def\mN{{\mathbb N}}

\def\mP{{\mathbb P}}

\def\mR{{\mathbb R}}

\def\sF{{\mathscr F}}

\def\sJ{{\mathscr J}}
\def\sK{{\mathscr K}}
\def\sL{{\mathscr L}}

\def\sQ{{\mathscr Q}}

\def\sU{{\mathscr U}}

\def\geq{\geqslant}
\def\leq{\leqslant}

\allowdisplaybreaks

\begin{document}

\title{Strong and weak convergence in the averaging principle for SDEs with H\"older coefficients}

\date{}

\author{Michael R\"ockner,\,\, Xiaobin Sun\,\, and\,\, Longjie Xie}

\address{Michael R\"ockner:
	Fakult\"{a}t f\"{u}r Mathematik, Universit\"{a}t Bielefeld, D-33501 Bielefeld, Germany, and Academy of Mathematics and Systems Science,
  Chinese Academy of Sciences (CAS), Beijing, 100190, P.R.China\\
	Email: roeckner@math.uni-bielefeld.de
}

\address{Xiaobin Sun:
	School of Mathematics and Statistics, Jiangsu Normal University,
	Xuzhou, Jiangsu 221000, P.R.China\\
	Email: xbsun@jsnu.edu.cn
}

\address{Longjie Xie:
	School of Mathematics and Statistics, Jiangsu Normal University,
	Xuzhou, Jiangsu 221000, P.R.China\\
	Email: longjiexie@jsnu.edu.cn
}

\thanks{
This work is supported in part by NSFC (No.11601196, 11701233, 11771187), NSF of Jiangsu (BK20170226) and the PAPD Project of Jiangsu Higher Education Institutions. Financial support of the DFG through CRC 1283 is gratefully acknowledged.}

\begin{abstract}
Using Zvonkin's transform and the Poisson equation in $\mR^d$ with a parameter,  we prove the averaging principle for stochastic differential equations with time-dependent H\"older continuous coefficients. Sharp convergence rates  with order  $(\alpha\wedge1)/2$ in the strong sense
 and $(\alpha/2)\wedge1$ in the weak sense are obtained, considerably extending the existing results in the literature. Moreover, we prove that the convergence of the multi-scale system  to the effective equation depends only on the regularity of the coefficients of the equation for the slow variable, and does not depend on the regularity of the coefficients of the equation for the fast component.

\bigskip


\end{abstract}

\maketitle

\section{Introduction}

In this paper, we consider the following stochastic slow-fast system in $\mR^{d_1}\times \mR^{d_2}$:
\begin{equation} \label{sde0}
\left\{ \begin{aligned}
&\dif X^{\eps}_t =\eps^{-1}b(X^{\eps}_t,Y^{\eps}_t)\dif t+\eps^{-1/2}\sigma(X^{\eps}_t,Y_t^\eps)\dif W^{1}_t,\qquad X^{\eps}_0=x\in\mR^{d_1},\\
&\dif Y^{\eps}_t =F(t, X^{\eps}_t, Y^{\eps}_t)\dif t+G(t,X_t^\eps,Y_t^\eps)\dif W^{2}_t,\,\,\,\quad\quad\quad Y^{\eps}_0=y\in\mR^{d_2},
\end{aligned} \right.
\end{equation}
where $d_1,d_2\geq 1$, $W^1_t$ and $W^2_t$ are $d_1$, $d_2$-dimensional independent standard Brownian motions both defined on some probability space $(\Omega,\sF,\mP)$, $b: \mR^{d_1}\times\mR^{d_2}\to\mR^{d_1}$, $F: \mR_+\times\mR^{d_1}\times\mR^{d_2}\to\mR^{d_2}$,  $\sigma: \mR^{d_1}\times\mR^{d_2}\to\mR^{d_1}\otimes\mR^{d_1}$ and $G: \mR_+\times\mR^{d_1}\times\mR^{d_2}\to\mR^{d_2}\otimes\mR^{d_2}$ are measurable functions, and the parameter $\eps>0$ represents the ratio between the  timescales of $X_t^\eps$ and $Y_t^\eps$ variables.
Such multiscale model  appears naturally in the theory of nonlinear oscillations, chemical kinetics, biology, climate dynamics and many other areas leading to a mathematical description involving `slow' and `fast' phase variables, see e.g. \cite{BR,  HKW, Ku, PS} and the references therein.
Usually, the underlying system (\ref{sde0}) is  difficult to deal with due to the two widely separated timescales and the cross interactions
of slow and fast modes. Hence, the asymptotic study of the behavior of the system  as $\eps\to0$ is of great interest and has attracted  much attentions in the past decades.

\vspace{1mm}
It is known that under suitable regularity assumptions on the coefficients, the slow part $Y_t^\eps$ will converge to the solution of the following reduced equation in $\mR^{d_2}$:
\begin{align}\label{sde1}
\dif \bar Y_t=\bar F(t,\bar Y_t)\dif t+\bar G(t,\bar Y_t)\dif W^2_t,\quad\bar Y_0=y,
\end{align}
where the new averaged coefficients are given by
\begin{align}\label{bf}
\bar F(t,y):=\!\int_{\mR^{d_1}}\!F(t,x,y)\mu^y(\dif x)\,\,\,\,\text{and}\,\,\,\,\bar G(t,y):=\!\sqrt{\int_{\mR^{d_1}}\!G(t,x,y)G(t,x,y)^{*}\mu^y(\dif x)}.
\end{align}
Here $G^{*}$ is the transpose of the matrix $G$, and $\mu^y(\dif x)$ is the unique invariant measure of the transition semigroup of the process
 $X_t^y$, which is the solution of the following frozen equation:
\begin{align}\label{sde2}
\dif X_t^y=b(X_t^y,y)\dif t+\sigma(X_t^y,y)\dif W_t^1,\quad X_0^y=x.
\end{align}
The effective dynamic (\ref{sde1}) then captures the evolution of the system (\ref{sde0}) over a long timescale, which does not depend on the fast variable any more and thus is much simpler than SDE (\ref{sde0}). This theory, known as the averaging principle, was first developed for  deterministic ordinary differential equations (ODEs for short) by Bogolyubov and Krylov \cite{KB}, and extended to the stochastic differential equations (SDEs for short) by Khasminskii \cite{K1}. We refer the readers to the book of Freidlin and Wentzell \cite{FW} for a comprehensive overview.

\vspace{1mm}
As a rule, the averaging method requires certain smoothness on both the original and the averaged coefficients. Various assumptions have been studied in order to guarantee the above convergence. Note that in the stochastic case, the convergence can be analyzed in two different ways: the strong convergence which provides pathwise asymptotic information for the system, and the weak convergence which gives convergence
the laws of the processes.
To the best of our knowledge, most of the results in the literature, both for the deterministic case and for the stochastic case, require at least local Lipschitz conditions on all the coefficients of system (\ref{sde0}), see e.g. \cite{GKK,Kh,Ki, Ki2, Li}. There is only one paper by Veretennikov \cite{V0} where weak convergence for the time-independent system (\ref{sde0}) was established under the assumptions that the drift coefficient $F$ in the slow equation is bounded and measurable with respect to the $y$ variable, and all the other coefficients are globally Lipschitz continuous.
%Thus, the averaging principle for SDEs with irregular coefficients has not been studied much.
Therefore, it seems that there are no studies of the averaging principle for SDEs which concentrates on H\"older coefficients.

\vspace{1mm}
On the other hand, in the papers mentioned above, no order of convergence in terms of  $\eps\to0$ is provided. But
for numerical purposes, it is important to know the rate of convergence of the slow variable to the effective dynamics. The main motivation comes from the well-known Heterogeneous Multi-scale Methods  used to approximate the slow component in system (\ref{sde0}), see e.g. \cite{Br3, ELV}. Moreover, the rate of convergence is also known to be very important for functional limit theorems in probability theory and homogenization, see e.g. \cite{KY,P-V,P-V2,WR}.
In this direction, the
strong convergence with order 1/2 and  weak convergence with order 1 are known to be optimal, see  \cite{G, KY2, L1,V, ZFWL}. As far as we know, all the known results in the literature concerning the rate of convergence require essentially at least $C^2_b$-regularity for all the coefficients, and none of them considered the fully coupled cases, i.e., the diffusion coefficient in the slow equation can not depend on the fast term.
We also mention that the averaging principle for stochastic partial differential equations and rates of convergence have also been widely studied, we refer to \cite{Br2, Br1,Ce,CF, DSXZ, FD} and the references therein.

\vspace{1mm}
The main aim of this work  is to develop a very general, robust and unified method for establishing the averaging principle, involving both strong and weak convergence, for the multi-scale system (\ref{sde0}) with {\bf irregular} coefficients, which leads to simplifications and extensions of the existing results. Unlike most previous publications, we mainly focus on the ``impact of noises" on the averaging principle for system (\ref{sde0}).
 More precisely,
we shall prove that under the non-degeneracy of the noises, the averaging principle holds for system (\ref{sde0}) with only H\"older continuous coefficients, see {\bf Theorem \ref{main4}}. Note that the deterministic system can even be ill-posed under such weak conditions on the coefficients. Moreover, we obtain the strong convergence rate with order $(\alpha\wedge1)/2$ and the weak convergence rate in the fully coupled case with order
$(\alpha/2)\wedge1$, where $\alpha>0$ is the H\"older index of the coefficients with respect to the slow component ($y$-variable), see {\bf Theorem \ref{main1}} and {\bf Theorem \ref{main2}} respectively.  In particular,
the convergence rates do not depend on the regularity of the coefficients with respect to the fast term ($x$-variable), which appear to be a new observation and which we think provides some new insight for understanding the averaging principle. See {\bf Remark \ref{br2}} and {\bf Remark \ref{br1}} for more detailed comparisons of our results with the previous publications on the subject.

\vspace{1mm}
The averaging principle for system (\ref{sde0}) is also known to be closely related to
the behavior of solutions for second-order parabolic and elliptic partial differential
equations, see \cite{Fr,KY,V0} and the references therein. In fact, the infinitesimal operator corresponding to $(X_t^\eps, Y_t^\eps)$ has the form
$$
\sL^{\eps}:=\eps^{-1}\sL_0(x,y)+\sL_1(x,y),
$$
where
\begin{align}
&\sL_0:=\sL_0(x,y):=\sum_{i,j} a_{ij}(x,y)\frac{\p^2}{\p x_i\p x_j}+b(x,y)\cdot\nabla_x,\label{l0}\\
&\sL_1:=\sL_1(x,y):=\sum_{i,j} H_{ij}(t,x,y)\frac{\p^2}{\p y_i\p y_j}+F(t,x,y)\cdot\nabla_y\label{l1}
\end{align}
with $a(x,y):=\sigma(x,y)\sigma(x,y)^*/2$ and $H(t,x,y):=G(t,x,y)G(t,x,y)^*/2$.  Given a $T>0$, consider the following Cauchy problem in $[0,T]\times\mR^{d_1}\times\mR^{d_2}$:
\begin{equation}
\left\{\begin{aligned}
&\p_tu^\eps(t,x,y)+\sL^{\eps} u^\eps(t,x,y)=\psi(y),\quad 0\leq t< T,\\
&u^\eps(T,x,y)=\varphi(y).\label{pde22}
\end{aligned}
\right.
\end{equation}
Using Theorem \ref{main2}, we can study the behavior of the solution $u^\eps$ to equation (\ref{pde22}) as $\eps\to0$. More precisely, we shall prove that $u^\eps(t,x,y)$ converges to the solution $\bar u(t,y)$ of the following reduced Cauchy problem in $[0,T]\times\mR^{d_2}$:
\begin{equation}\left\{\begin{array}{l}\label{pde11}
\displaystyle
\p_t\bar u(t,y)+\bar{\mathscr{L}}\bar u(t,y)=\psi(y),\quad 0\leq t< T,\\
\bar u(T,y)=\varphi(y),
\end{array}\right.
\end{equation}
where $\psi$ is a bounded measurable function, $\varphi$ is bounded continuous,  and $\bar{\mathscr{L}}$ is the infinitesimal generator of the effective SDE (\ref{sde1}), i.e.,
\begin{align}\label{lby}
\bar\sL:=\sum_{i,j}\bar H_{ij}(t,y)\frac{\p^2}{\p y_i\p y_j}+\bar F(t,y)\cdot\nabla_y
\end{align}
with $\bar H(t,y):=\bar G(t,y)\bar G(t,y)^*/2$, and $\bar F, \bar G$ are as defined in (\ref{bf}). The main result in this direction is given by {\bf Theorem \ref{main3}}.

\vspace{1mm}
As mentioned before, the argument that we shall use is rather simple insofar as it does not involve the classical time discretisation procedure, which is commonly used in the literature to prove the averaging principle. Two ingredients are crucial in our proof: Zvonkin's transformation  and the Poisson equation in the whole space. First of all, due to the low regularity of the coefficients, we shall use Zvonkin's argument to transform the equation for $Y_t^\eps$ and $\bar Y_t$ into new ones. Such technique was first developed in \cite{Zv} and is now widely used to study the strong well-posedness for SDEs with singular coefficients, see e.g. \cite{Kr-Ro,XZ,XZ2,Zh1}.
Then we use the  Poisson equation with a parameter to prove both the strong and weak convergence for system (\ref{sde0}). Here we adopt and improve the idea used in \cite{Br1}, where  the convergence rate in the averaging principle for SPDEs with smooth coefficients with the fast equation not depending on the slow component was studied.
More precisely, we shall study the following Poisson equation in $\mR^{d_1}$:
\begin{align}
\sL_0(x,y)u(x,y)=-f(x,y),\quad x\in\mR^{d_1},  \label{pde1}
\end{align}
where $y\in\mR^{d_2}$ is a parameter and $\sL_0(x,y)$ is defined by (\ref{l0}).
We note that there is no boundary condition. When the equation is formulated in a compact set, the corresponding theory is well known. However, equation (\ref{pde1}) in the whole space $\mR^{d_1}$ has been studied only very recently, and it turns out to be  very useful in the theory of the averaging principle, diffusion approximation and other limit theorems, see the series of papers \cite{BSV,P-V,P-V2,Ve}. We shall derive estimates for the solution of (\ref{pde1}) in terms of explicit conditions on the coefficients as well as the right hand side, see {\bf Theorem \ref{popde}}, which generalizes the results in  \cite{P-V,P-V2} and is of independent interest.

\vspace{1mm}
The paper is organized as follows. In Section 2, we state our main results. Section 3 is devoted to the study of the Poisson equation in the whole space with a parameter. The proofs  of strong convergence and weak convergence are given in Section 4 and Section 5, respectively.

\vspace{2mm}
To end this section, we introduce some notations. Let $\mN:=\{0,1,\cdots\}$ and $\mN^*:=\{1,2\cdots\}$. For $\beta\in (0,1)$, let $C^{\beta}(\mR^d)$ be the usual local H\"older space.  For $\beta\in\mN^*$, without abuse of notation, we denote by $C^{\beta}(\mR^d)$ the space of all functions $f$ whose $\beta-1$ order derivative $\partial^{\beta-1}f$ is Lipschitz continuous.  While when  $\beta\in(0,\infty)\setminus\mN^*$, $C^{\beta}(\mR^d)$ consists of all functions satisfying $f\in C^{[\beta]}(\mR^d)$ and $\partial^{[\beta]}f\in C^{\beta-[\beta]}(\mR^d)$, where $[\beta]$ denotes the largest integer which is smaller than $\beta$. For $\beta>0$, we denote by $C^{\beta}_b(\mR^d)$ the space of all functions $f\in C^{\beta}(\mR^d)$ whose $i$-order derivative $\partial^{i}f$ is bounded for any $0\leq i\leq ([\beta]-1)\vee0$.

Given a function $f$ and $\gamma_1,\gamma_2,\gamma_3\in(0,\infty)$, we shall consider the following three  cases: \\
(i) $f$ is defined on $\mR^{d_1+d_2}$, i.e., $f$ is a function with variable $x$ and $y$: we write $f\in C_b^{\gamma_1,\gamma_2}$ if $f\in C_b^{\gamma_1}(\mR^{d_1};C_b^{\gamma_2}(\mR^{d_2}))$, and $f\in C_{loc,y}^{\gamma_1,\gamma_2}$ means $f\in C_b^{\gamma_1}(\mR^{d_1};C^{\gamma_2}(\mR^{d_2}))$;\\
(ii) $f$ is defined on $\mR_+\times\mR^{d_1+d_2}$, i.e., $f$ is a function of $t,x$ and $y$: we write $f\in C_b^{\gamma_3,\gamma_1,\gamma_2}$ if for every $t>0$, $f(t,\cdot,\cdot)\in C_b^{\gamma_1,\gamma_2}$ and for every $(x,y)$, $f(\cdot,x,y)\in C_b^{\gamma_3}(\mR_+)$; similarly, $f\in C_{loc,y}^{\gamma_3,\gamma_1,\gamma_2}$ means that $f(t,\cdot,\cdot)\in C_{loc,y}^{\gamma_1,\gamma_2}$ and $f(\cdot,x,y)\in C_b^{\gamma_3}(\mR_+)$;\\
(iii) $f$ is defined on $\mR_+\times\mR^{d_2}$, i.e.,  $f$ is a function of $t$ and $y$:  we write $f\in C_b^{\gamma_3,\gamma_2}$ if for every $t>0$, $f(t,\cdot)\in C_b^{\gamma_2}$ and for every $y\in\mR^{d_2}$, $f(\cdot,y)\in C_b^{\gamma_3}(\mR_+)$.

\section{Assumptions and main results}

Let us first introduce some basic assumptions. We shall assume the following non-degeneracy conditions on the diffusion coefficients:

\vspace{2mm}
\noindent{\bf (H$\sigma$):} The coefficient $a=\sigma\sigma^*$ is non-degenerate in $x$ uniformly with respect to $y$, i.e., there exists a $\lambda>1$ such that  for any $x\in\mR^{d_1}$ and $y\in\mR^{d_2}$,
$$
\lambda^{-1}|\xi|^2\leq a^{ij}(x,y)\xi_i\xi_j\leq \lambda|\xi|^2,\ \ \forall\xi\in\mR^{d_1}.
$$

\vspace{2mm}
\noindent{\bf (H$_{G}$):} The coefficient $H=GG^*$ is non-degenerate in $y$ uniformly with respect to $(t,x)$, i.e., there exists a $\lambda>1$ such that for any $(t,x)\in\mR_+\times\mR^{d_1}$  and $y\in\mR^{d_2}$,
$$
\lambda^{-1}|\xi|^2\leq H^{ij}(t,x,y)\xi_i\xi_j\leq\lambda|\xi|^2,\ \ \forall\xi\in\mR^{d_2}.
$$

For the existence of an invariant measure for the frozen SDE (\ref{sde2}), we assume the following very weak recurrence condition (see \cite{P-V2, Ve5}):

\vspace{2mm}
\noindent{\bf (H$b$):}\qquad\qquad\qquad\quad $\lim_{|x|\to\infty}\sup_y \<x,b(x,y)\>=-\infty$.

\vspace{2mm}
Below, we state our main results concerning the strong and weak convergence for the  averaging principle for system (\ref{sde0}) separately.

\subsection{Strong convergence}

The following is the first main result of this paper.

\bt\label{main1}
Let {\bf (H$\sigma$)}-{\bf (H$_{G}$)}-{\bf (H$b$)} hold,  and let
\begin{align}\label{gy}
G(t,x,y)\equiv G(t,y).
\end{align}
Assume that $\sigma\in C_b^{1,1}$, $b\in C_{b}^{\delta,\alpha}$ and  $F\in C_{b}^{\alpha/2,\delta,\alpha}$, $G\in C_b^{\alpha/2,1}$ with $0<\delta,\alpha\leq1$. Then we have for any $T>0$,
\begin{align}\label{stro}
\sup_{t\in[0,T]}\mE|Y_t^\eps-\bar Y_t|^2\leq C_T\,\eps^{\alpha\wedge1},
\end{align}
where $C_T>0$ is a constant independent of $\delta$.
\et

We point out that under our assumptions, the strong well-posedness for system (\ref{sde0}) was obtained by \cite{Ver} or \cite[Theorem 1.3]{Zh1}, and the invariant measure $\mu^y(\dif x)$ for SDE (\ref{sde2}) exists and is unique, see \cite[Theorem 1.2]{XZ} or \cite[Theorem 2.9]{XZ2}. Meanwhile, we shall show that the averaged drift $\bar F$ defined in (\ref{bf}) is also H\"older continuous, i.e., $\bar F\in C_b^{\alpha/2,\alpha}$ (see Lemma \ref{hff} below). Thus, there exists a unique strong solution $\bar Y_t$ to SDE (\ref{sde1}).

 Let us list some important comments  to explain our result.

\br\label{br2}
We first point out that the independence of $G$ with respect to the $x$-variable in assumption (\ref{gy}) is necessary. Otherwise, the strong convergence for SDE (\ref{sde0}) may not be true, cf. \cite{L1,V0}.

\vspace{1mm}
\noindent{\bf (1) [Singular coefficients].} We do not make any Lipschitz-type assumptions on the drift coefficients $b$ and $F$. This is due to the regularization effect of the non-degenerate noises. Note that if $\sigma=0$ or $G=0$, the system (\ref{sde0}) may even be ill-posed with only H\"older coefficients.

\vspace{1mm}
\noindent{\bf(2) [Sharp order].}  Taking $\alpha=1$ in (\ref{stro}), we can obtain the strong convergence with order $1/2$. Thus we get the optimal rate under much weaker regularity conditions both on the diffusion and the drift coefficients than the known results in the literature.  Meanwhile, when $0<\alpha<1$, we also get that the averaging principle holds with a strong convergence rate $\alpha/2$, which to the best of our knowledge is  new. Moreover, we allow the coefficients to be time-dependent, which appears to have not been studied before in estimating the rate of convergence.

\vspace{1mm}
\noindent{\bf(3) [Dependence of convergence].} Note that the convergence rate $(\alpha\wedge1)/2$ does not depend on the index $\delta$. This suggests that the convergence in the averaging principle replies only on the regularity of the coefficients  with respect to the $y$ (slow) variable, and does not depend on the regularity  with respect to the $x$ (fast) variable, which we think provides some new insight for understanding the averaging principle.
\er

By a localization technique as in \cite[Corollary 2.6]{XZ2}, we can drop the boundness condition on the coefficients with respect to the slow variable.
\bt\label{main4}
Let {\bf (H$\sigma$)}-{\bf (H$_{G}$)}-{\bf (H$b$)} hold, and let
\begin{align*}
G(t,x,y)\equiv G(t,y).
\end{align*}
Assume $\sigma\in C^{1,1}_{loc,y}$, $G\in C^{\alpha/2,1}_{loc,y}$ and $b\in C^{\delta,\alpha}_{loc,y}$, $F\in C^{\alpha/2,\delta,\alpha}_{loc,y}$ with $0<\delta,\alpha\leq1$, and that the following moment estimate holds:

\begin{enumerate}
	\item [{\bf (H$^{M}$)}] For any $T>0$, there exists a $\beta>2$ such that
	\begin{eqnarray*}
		\sup_{\eps\in(0,\eps_0)}\mE\Big[\sup_{t\in[0,T]}(|Y_t^\eps|^\beta+|\bar Y_t|^\beta)\Big]\leq C<\infty,
	\end{eqnarray*}
	where $\eps_0>0$ and $C>0$ is a constant depending on $T, |x|,|y|$, where $x,y$ are the initial conditions in (\ref{sde0}).
\end{enumerate}
 Then we have for any $T>0$,
$$
\lim_{\eps\to 0}\sup_{t\in[0,T]}\mE|Y_t^\eps-\bar Y_t|^2=0.
$$
\et
\br
Local conditions imposed on the coefficients allow functions to have certain growth at infinity.
The advantage of Theorem \ref{main4} lies in that, we only need to show the a priori moment estimate {\bf (H$^{M}$)} in order to guarantee  the strong convergence in the averaging principle for SDE (\ref{sde0}) with only local H\"older continuous drifts.
\er

\subsection{Weak convergence}
In the above results, the assumptions $\sigma\in C_b^{1,1}$ and $G\in C_b^{\alpha/2,1}$ are mainly needed to ensure the strong well-posedness for system (\ref{sde0}). Now, we state our main result concerning the weak convergence of system (\ref{sde0}) under weaker conditions on the diffusion coefficients.

\bt[Weak convergence]\label{main2}
Let {\bf (H$\sigma$)}-{\bf (H$_{G}$)}-{\bf (H$b$)} hold true.
Assume that $\sigma, b\in C_b^{\delta,\alpha}$ and  $F, G\in C_{b}^{\alpha/2,\delta,\alpha}$ with $0<\delta,\alpha\leq2$. Then for any $T>0$
and every $\varphi\in C_b^{2+\alpha}(\mR^{d_2})$, we have
\begin{align}\label{weak}
\sup_{t\in[0,T]}\Big|\mE[\varphi(Y_t^\eps)]-\mE[\varphi(Y_t)]\Big|\leq C_T\,\eps^{(\alpha/2)\wedge1},
\end{align}
where $C_T>0$ is a constant independent of $\delta$.
\et

We now give some comments to explain the above result.

\br\label{br1}
Note that here the diffusion coefficient $G$ in the slow equation can also depend on the fast variable $x$.

\vspace{1mm}
\noindent{\bf (1) [Singular coefficients].}
Due to the non-degeneracy of the noises, it is well-known that the system (\ref{sde0}) is weakly well-posed under our conditions.
%%%This is the main reason why we can assume weaker conditions on the diffusion coefficients to study the weak convergence.
This is the main reason why we can assume weaker conditions on the diffusion coefficients to prove the above weak convergence.

\vspace{1mm}
\noindent{\bf(2) [Sharp order].}  Taking $\alpha=2$ in (\ref{weak}), we obtain  the optimal weak convergence rate $1$.
Our result generalizes the known results in the literature by allowing the coefficients to be time-dependent, and more importantly, to be fully coupled, i.e.,  the diffusion coefficient in the slow equation can depend on the fast variable, which appears to have not been considered before in estimating the rate of convergence.
Meanwhile, when $0<\alpha<2$, we also get that the weak averaging principle holds with convergence rate $\alpha/2$, which also appears to be new.

\vspace{1mm}
\noindent{\bf(3) [Dependence of convergence].} As before, the weak convergence relies only on the regularity of all the coefficients  with respect to the $y$ (slow) variable, since the rate  $(\alpha/2)\wedge1$ does not depend on the index $\delta$.
\er

As a direct consequence of Theorem \ref{main2}, we have the following result concerning the limit behavior of parabolic equations.

\bt\label{main3}
Suppose the assumptions in Theorem \ref{main2} hold, $\psi$ is bounded measurable and $\varphi$ is bounded continuous. Let $u^\eps$ be the solution to equation (\ref{pde22}). Then for every $t>0$, $x\in\mR^{d_1}$ and $y\in\mR^{d_2}$, the limit
$$
\lim_{\eps\to0}u^\eps(t,x,y)=:u(t,y)
$$
exists, and the function $u(t,y)$ is the unique solution of the Cauchy problem (\ref{pde11}).
\et

\section{Poisson equation in $\mR^{d_1}$ with a parameter}

This section is devoted to studying the Poisson equation (\ref{pde1}) in the whole space. We are looking for a solution $u$ for (\ref{pde1}) which grows at most polynomial in $x$ as $|x|\to\infty$, and the main problem addressed here is the regularity of the solution $u$ with respect to the parameter $y$.  Throughout this section, we shall always assume
{\bf (H$\sigma$)} and {\bf (H$b$)} to hold.
Let us point out  that there is no boundary condition. As a result, the solution turns out to be defined up to an additive constant, since $\sL(x,y)1\equiv0$. To fix this constant, it is convenient to make the following ``centering" assumption on the right-hand side:
\begin{align}\label{cen}
\int_{\mR^{d_1}}f(x,y)\mu^y(\dif x)=0,\quad\forall y\in\mR^{d_2},
\end{align}
which is analogous to the centering in the standard Central Limit Theorem, see \cite{P-V,P-V2} for more details.

We shall essentially use the strategy implemented in \cite{P-V2}, where the fundamental solution was used  to study the equation (\ref{pde1}).
More precisely, note that $\sL_0(x,y)$ can be viewed as the infinitesimal generator of the process $X_t^y(x)$, which is the unique strong solution for the frozen SDE (\ref{sde2}). As a result, the solution $u$ to equation (\ref{pde1}) should have the
following probabilistic representation:
\begin{align}\label{pro}
u(x,y)=\int_0^\infty\mE f\big(X_t^y(x),y\big)\dif t.
\end{align}
As we shall see below,
under our assumptions, $X_t^y(x)$ admits a density function $p_t(x,x';y)$, which is also the unique fundamental solution for the operator $\sL_0(x,y)$.
Let  $T_tf(x,y)$denotes  the semigroup corresponding to $X_t^y(x)$, i.e.,
$$
T_tf(x,y):=\mE \big(f(X_t^y(x),y\big)=\int_{\mR^{d_1}}p_t(x,x';y)f(x',y)\dif x'.
$$
Then we can write
\begin{align}\label{uu}
u(x,y)=\int_0^1T_tf(x,y)\dif t+\int_1^\infty T_tf(x,y)\dif t.
\end{align}
Thus, we need to study the behavior of $T_tf$ as well as its first and second order derivatives with respect to the $y$-variable both near $t=0$ and as $t\to\infty$.
The following is the main result of this section.

\bt\label{popde}
Let {\bf (H$\sigma$)} and {\bf (H$b$)} hold. Assume that $a, b\in C_b^{\delta,\ell}$ with $0<\delta\leq 1$ and $\ell=0,1,2$. Then for every function $f\in C_b^{\delta,\ell}$ satisfying (\ref{cen}), there exists a unique solution $u$ to (\ref{pde1}) such that for any $y\in\mR^{d_2}$, $u(\cdot,y)\in C^{2}$ and for any $x\in\mR^{d_1}$, $u(x,\cdot)\in C_b^{\ell}$. Moreover, there exists a constant $m>0$ such that for any $y\in\mR^{d_2}$,
\begin{align}\label{key1}
|u(x,y)|+|\nabla_xu(x,y)|+|\nabla_x^2u(x,y)|\leq C_0\|f\|_{C_b^{\delta,0}}(1+|x|^m),
\end{align}
and when $\ell=1$,
\begin{align}\label{key3}
|\nabla_yu(x,y)|\leq &C_0\Big[\big(\|a\|_{C_b^{\delta,0}}+\|b\|_{C_b^{\delta,0}}\big)\|f\|_{C_b^{\delta,1}}\no\\
&+\big(\|a\|_{C_b^{\delta,1}}+\|b\|_{C_b^{\delta,1}}\big)\|f\|_{C_b^{\delta,0}}\Big](1+|x|^m),
\end{align}
and when $\ell=2$,
\begin{align}\label{key2}
|\nabla^2_yu(x,y)|&\leq C_0\Big[\big(\|a\|_{C_b^{\delta,0}}+\|b\|_{C_b^{\delta,0}}\big)\|f\|_{C_b^{\delta,2}}+\big(\|a\|_{C_b^{\delta,1}}+\|b\|_{C_b^{\delta,1}}\big)\|f\|_{C_b^{\delta,1}}\no\\
&+\Big(\big(\|a\|_{C_b^{\delta,1}}+\|b\|_{C_b^{\delta,1}}\big)^2+\big(\|a\|_{C_b^{\delta,2}}+\|b\|_{C_b^{\delta,2}}\big)\Big)\|f\|_{C_b^{\delta,0}}\Big](1+|x|^m),
\end{align}
where $C_0$ is a positive constant depending only on $\lambda,d_1,d_2$ and $\|a\|_{C_b^{\delta,0}}, \|b\|_{C_b^{\delta,0}}$.
\et

\br
Concerning estimates (\ref{key3})-(\ref{key2}), usually one does not care about the dependence of  constants  on the right hand side with respect to the norms of the coefficients. But this will be very important below for us to get the sharp rate of convergence for system (\ref{sde0}) with only H\"older continuous coefficients. More precisely, since we assume the coefficients  belong to the space $C_b^{\delta,\alpha}$ in Theorem \ref{main1} and Theorem \ref{main2}, we need to keep track of the dependence of the constant on the right hand side of (\ref{key3})-(\ref{key2}) with respect to the higher order norms of the coefficients $a,b$ as well as of the potential term $f$.
\er
\begin{proof}[Proof of Theorem \ref{popde}]
	The existence and uniqueness of the solution $u$ to (\ref{pde1}) are well-known under the above conditions. Meanwhile,  by regarding $y\in\mR^{d_2}$ as a parameter, the estimate (\ref{key1}) is true since all coefficients are bounded uniformly in the $y$ variable, see e.g. \cite{P-V,P-V2}. Concerning estimate (\ref{key3}), we have by (\ref{uu}) and Lemma \ref{dyf} below that for any $k\in\mR_+$, there exist constants $C_1,m>0$ such that
	\begin{align*}
	|\nabla_y u(x,y)|&\leq \int_0^1|\nabla_y T_tf(x,y)|\dif t+\int_1^\infty |\nabla_y T_tf(x,y)|\dif t\\
	&\leq C_1 \Big[\big(\|a\|_{C_b^{\delta,0}}+\|b\|_{C_b^{\delta,0}}\big)\|f\|_{C_b^{\delta,1}}+\big(\|a\|_{C_b^{\delta,1}}+\|b\|_{C_b^{\delta,1}}\big)\|f\|_{C_b^{\delta,0}}\Big]\\
&\quad\times\left(1+\int_1^\infty\frac{(1+|x|^m)}{(1+t)^k}\dif t\right),
	\end{align*}
which in turn yields the desired result. The estimate (\ref{key2}) can be proved similarly. The proof is finished.
\end{proof}

Below, we proceed to study the first and second order derivatives of $T_tf(x,y)$ with respect to  the $y$-variable in the following two subsections. We provide the explicit dependence on all higher order norms of the coefficients involved.

\subsection{First order derivative with respect to $y$} Let us first recall some classical results concerning the fundamental solution $p_t(x,x';y)$, see \cite[Theorem 2.3]{CHXZ}, \cite[Chapter IV]{La-So-Ur} and \cite[Proposition 3]{P-V2}.

\bl
Assume {\bf (H$\sigma$)} holds and let $T>0$. Let $a,b\in C_b^{\delta,0}$ with  $0<\delta\leq 1$. Then
for every $\ell=0,1,2$ and any $0<t\leq T$, we have
\begin{align}\label{ppp1}
|\nabla_x^\ell p_t(x,x';y)|\leq C_{T}t^{-(d+\ell)/2}\exp\big(-c_0|x-x'|^2/t\big),
\end{align}
and for every $x_1,x_2\in\mR^{d_1}$ and $0<\delta'\leq \delta$,
\begin{align}\label{ppp2}
|\nabla_x^2 p_t(x_1,x';y)&-\nabla_x^2 p_t(x_2,x';y)|\leq C_{T}|x_1-x_2|^{\delta'}t^{-(d+2+\delta')/2}\no\\
&\quad\times\Big(\exp\big(-c_0|x_1-x'|^2/t\big)+\exp\big(-c_0|x_2-x'|^2/t\big)\Big),
\end{align}
where $C_T,c_0>0$ are constants independent of $y$.

If we further assume {\bf (H$b$)} holds, then for any  $k, j\in\mR_+$, there exists a constant $m>0$ such that for all $t\geq 1$, $x,x'\in\mR^{d_1}$ and $y\in\mR^{d_2}$,
\begin{align}\label{pp1}
|p_t(x,x';y)|\leq C_1\frac{1+|x|^m}{(1+|x'|^{j})},
\end{align}
and for $\ell=1,2$,
\begin{align}\label{pp2}
|\nabla^\ell_xp_t(x,x';y)|\leq C_2\frac{1+|x|^m}{(1+t)^k(1+|x'|^j)}.
\end{align}
Moreover, the limit
$$
p_\infty(x',y):=\lim_{t\to\infty}p_t(x,x';y)
$$
exists and is independent of $x$, and for every $k, j\in\mR_+$, there exists a constant $m>0$ such that for any $y\in\mR^{d_2}$,
\begin{align}\label{pin1}
|p_\infty(x',y)|\leq \frac{C_3}{1+|x'|^j},
\end{align}
and
\begin{align}\label{pin2}
|p_t(x,x';y)-p_\infty(x',y)|\leq C_4\frac{1+|x|^m}{(1+t)^k(1+|x'|^j)}.
\end{align}
The above positive constants $C_i (i=1,\cdots,4)$ depend only on $\lambda,d_1,d_2$ and $\|a\|_{C_b^{\delta,0}}, \|b\|_{C_b^{\delta,0}}$.
\el

To study the regularity of $T_tf$ with respect to the $y$-variable, we first consider the case where $f(x,y)\equiv g(x)$, i.e., the function $f$ does not depend on the parameter $y$, and
\begin{align}\label{cecece}
\int_{\mR^d}g(x)\mu^y(\dif x)\equiv0,\quad\forall y\in\mR^{d_2}.
\end{align}
To shorten the notation, we write for $\ell=1,2$,
$$
\frac{\p^\ell\sL_0}{\p y^\ell}(x,y):=\sum_{i,j} \p_y^\ell a_{ij}(x,y)\frac{\p^2}{\p x_i\p x_j}+\p_y^\ell b(x,y)\cdot\nabla_x.
$$
We have the following result.

\bl
Let {\bf (H$\sigma$)}-{\bf (H$b$)} and (\ref{cen}) hold. Assume that $a, b\in C_b^{\delta,1}$ and $g\in C_b^\delta$ with $0<\delta\leq 1$. Then we have
\begin{align}\label{d1}
\nabla_yT_tg(x,y)=\int_0^t\!\!\int_{\mR^{d_1}}p_{t-s}(x,z;y)\frac{\p\sL_0}{\p y}(z,y)T_sg(z,y)\dif z\dif s.
\end{align}
Moreover, for any $0<t\leq 2$,
\begin{align}\label{es1}
|\nabla_yT_tg(x,y)|\leq C_0\Big(\|a\|_{C_b^{\delta,1}}+\|b\|_{C_b^{\delta,1}}\Big)\|g\|_{C_b^\delta},
\end{align}
and for any $k\in\mR_+$ , there exists a constant $m>0$ such that for all $t>2$,
\begin{align}\label{es2}
|\nabla_yT_tg(x,y)|\leq C_0\Big(\|a\|_{C_b^{\delta,1}}+\|b\|_{C_b^{\delta,1}}\Big)\|g\|_{C_b^\delta} \frac{(1+|x|^m)}{(1+t)^k},
\end{align}
where $C_0>0$ is a constant depending only on $\lambda,d_1,d_2$ and $\|a\|_{C_b^{\delta,0}}, \|b\|_{C_b^{\delta,0}}$.
\el
\begin{proof}
	The equality (\ref{d1}) has been proved in \cite[Theorem 10]{P-V2} under sightly stronger assumptions on the coefficients. Let us show that the right hand side is indeed well-defined under our conditions.
	In fact, since $a, b\in C_b^{\delta,1}$, the operator $\p\sL_0/\p y$ is meaningful.
	On the other hand, since $g\in C_b^\delta$, we can derive by (\ref{ppp1}) that for $\ell=1,2$ and any $0<s\leq 2$,
	\begin{align}\label{can}
	\nabla^\ell_zT_sg(z,y)&=\int_{\mR^{d_1}}\nabla_z^\ell p_s(z,x';y)\big[g(x')-g(z)\big]\dif x'\no\\
	&\leq C_1\|g\|_{C_b^\delta}\int_{\mR^{d_1}}s^{-(d+\ell)/2}\exp\big(-c_0|z-x'|^2/s\big)\cdot|x'-z|^\delta\dif x'\no\\
	&\leq C_1\|g\|_{C_b^\delta}s^{(\delta-\ell)/2},
	\end{align}
	and for any $s>1$, we have by (\ref{pp2}) that
	\begin{align}\label{can2}
	\nabla^\ell_zT_sg(z,y)\leq C_1\int_{\mR^{d_1}}\frac{1+|z|^m}{(1+s)^k(1+|x'|^j)}|g(x')|\dif x'\leq C_1\|g\|_{C_b^\delta}\frac{1+|z|^m}{(1+s)^k},
	\end{align}
	where $C_1>0$ is a constant independent of $s$ and $y$.
As a result,
$$
\frac{\p\sL_0}{\p y}(z,y)T_sg(z,y)
$$
makes sense, and estimate (\ref{es1}) follows directly. Below, we proceed to show (\ref{es2}).
As a consequence of (\ref{cecece}), we have (see \cite[(28)]{P-V2})
\begin{align*}
\nabla_yT_tg(x,y)&=\int_0^t\!\!\int_{\mR^{d_1}}p_{t-s}(x,z;y)\frac{\p\sL_0}{\p y}(z,y)T_sg(z,y)\dif z\dif s\\
&\quad-\int_0^\infty\!\!\int_{\mR^{d_1}}p_\infty(z;y)\frac{\p\sL_0}{\p y}(z,y)T_sg(z,y)\dif z\dif s.
\end{align*}
We further write
\begin{align*}
\nabla_yT_tg(x,y)&= \int_0^{t/2}\!\!\int_{\mR^{d_1}}\big[p_{t-s}(x,z;y)-p_\infty(z;y)\big]\frac{\p\sL_0}{\p y}(z,y)T_sg(z,y)\dif z\dif s\\
&\quad+\int_{t/2}^t\int_{\mR^{d_1}}p_{t-s}(x,z;y)\frac{\p\sL_0}{\p y}(z,y)T_sg(z,y)\dif z\dif s\\
&\quad-\int_{t/2}^\infty\!\int_{\mR^{d_1}}p_\infty(z;y)\frac{\p\sL_0}{\p y}(z,y)T_sg(z,y)\dif z\dif s=:\cI_1+\cI_2+\cI_3.
\end{align*}
For the first term, we have by (\ref{pin2}) that for any $k\in\mR_{+}$,
\begin{align*}
\cI_1&\leq \left(\int_0^{1}+\int_1^{t/2}\right)\!\int_{\mR^{d_1}}C_2\frac{1+|x|^m}{(1+t-s)^k(1+|z|^j)}\Big|\frac{\p\sL_0}{\p y}(z,y)T_sg(z,y)\Big|\dif z\dif s.
\end{align*}
Using (\ref{can}) and (\ref{can2}) , we can derive that
\begin{align*}
\cI_1&\leq C_2\Big(\|a\|_{C_b^{\delta,1}}+\|b\|_{C_b^{\delta,1}}\Big)\|g\|_{C_b^\delta}\bigg(\int_0^{1}\!\!\int_{\mR^{d_1}}\frac{1+|x|^m}{(1+t-s)^k(1+|z|^j)}s^{\delta/2-1}\dif z\dif s\\
&\quad+\int_1^{t/2}\!\!\int_{\mR^{d_1}}\frac{1+|x|^m}{(1+t-s)^k(1+|z|^j)}\frac{1+|z|^m}{(1+s)^k}\dif z\dif s\bigg)\\
&\leq C_2\Big(\|a\|_{C_b^{\delta,1}}+\|b\|_{C_b^{\delta,1}}\Big)\|g\|_{C_b^\delta}\frac{(1+|x|^m)}{(1+t)^k},
\end{align*}
where we choose $j>d_1+m$ in the last inequality. As for the second term, using (\ref{pp1}) and (\ref{pp2}), we have
\begin{align*}
\cI_2&=\int^{t/2}_0\!\!\int_{\mR^{d_1}}p_{s}(x,z;y)\frac{\p\sL_0}{\p y}(z,y)T_{t-s}g(z,y)\dif z\dif s\\
&\leq C_3\Big(\|a\|_{C_b^{\delta,1}}+\|b\|_{C_b^{\delta,1}}\Big)\|g\|_{C_b^\delta}\bigg(\int^{1}_0\!\!\int_{\mR^{d_1}}s^{-d/2}\exp\big(-c_0|x-z|^2/s\big)\\
&\qquad\qquad\qquad\times\frac{1+|z|^m}{(1+(t-s))^k}\dif z\dif s+\int^{t/2}_1\!\!\int_{\mR^{d_1}}\frac{1+|x|^m}{(1+|z|^{j})}\frac{1+|z|^m}{(1+(t-s))^k}\dif z\dif s\bigg)\\
&\leq C_3\Big(\|a\|_{C_b^{\delta,1}}+\|b\|_{C_b^{\delta,1}}\Big)\|g\|_{C_b^\delta}\frac{(1+|x|^m)}{(1+t)^k}.
\end{align*}
Finally, we have by (\ref{pin1}) that
	\begin{align*}
	\cI_3&\leq C_4\Big(\|a\|_{C_b^{\delta,1}}+\|b\|_{C_b^{\delta,1}}\Big)\|g\|_{C_b^\delta}\int_{t/2}^\infty\!\int_{\mR^{d_1}}\frac{1}{1+|z|^j}\frac{1+|z|^m}{(1+s)^k}\dif z\dif s\\
	&\leq C_4\Big(\|a\|_{C_b^{\delta,1}}+\|b\|_{C_b^{\delta,1}}\Big)\|g\|_{C_b^\delta}\frac{1}{(1+t)^k}\leq C_4\Big(\|a\|_{C_b^{\delta,1}}+\|b\|_{C_b^{\delta,1}}\Big)\|g\|_{C_b^\delta}\frac{(1+|x|^m)}{(1+t)^k}.
	\end{align*}
The proof is finished.
\end{proof}

\subsection{Second order derivative with respect to $y$}

We shall need the following regularity result for $\nabla_yT_tg(x,y)$ with respect to $x$ to study the second order derivative of $T_tg$ with respect to the parameter $y$.

\bl
Let {\bf (H$\sigma$)}, {\bf (H$b$)} and (\ref{cen}) hold. Assume that $a, b\in C_b^{\delta,1}$ and $g\in C_b^\delta$ with  $0<\delta\leq 1$. Then  for every $y\in\mR^{d_2}$, we have
$\nabla_yT_tg(\cdot,y)\in C^2$.
Moreover, for any $0<t\leq 2$ and $\ell=1,2$,
\begin{align}\label{xy1}
|\nabla^\ell_x\nabla_yT_tg(x,y)|\leq C_0t^{(\delta-\ell)/2}\Big(\|a\|_{C_b^{\delta,1}}+\|b\|_{C_b^{\delta,1}}\Big)\|g\|_{C_b^\delta},
\end{align}
and for any $k\in\mR_+$, there exists a constant $m>0$ such that  for any $t>2$,
\begin{align}\label{xy2}
|\nabla^\ell_x\nabla_yT_tg(x,y)|\leq C_0\Big(\|a\|_{C_b^{\delta,1}}+\|b\|_{C_b^{\delta,1}}\Big)\|g\|_{C_b^\delta}\frac{(1+|x|^m)}{(1+t)^k},
\end{align}
where $C_0>0$ is a constant depending only on $\lambda,d_1,d_2$ and $\|a\|_{C_b^{\delta,0}}, \|b\|_{C_b^{\delta,0}}$.
\el
\begin{proof}
We only  prove the estimates (\ref{xy1}) and (\ref{xy2}) when  $\ell=2$, the case $\ell=1$ can be proved similarly and is easier since it involves less singularities.
Recall that
\begin{align*}
\nabla_yT_tg(x,y)= \int_0^{t}\!\!\int_{\mR^d}p_{t-s}(x,z;y)\frac{\p\sL_0}{\p y}(z,y)T_sg(z,y)\dif z\dif s.
\end{align*}
By the H\"older assumption on the coefficients and (\ref{ppp2}), it is easy to check that the function
$$
z\to\frac{\p\sL_0}{\p y}(\cdot,y)T_sg(\cdot,y)
$$
is $\delta'$-H\"older continuous for any $\delta'<\delta$, i.e., for any $x,z\in\mR^{d_1}$ and $s\leq 2$, there exists a constant $C_1>0$ such that
\begin{align*}
&\bigg|\frac{\p\sL_0}{\p y}(x,y)T_sg(x,y)-\frac{\p\sL_0}{\p y}(z,y)T_sg(z,y)\bigg|\\
&\leq C_1\Big(\|a\|_{C_b^{\delta,1}}+\|b\|_{C_b^{\delta,1}}\Big)\|g\|_{C_b^\delta}|x-z|^{\delta'}s^{\frac{\delta-\delta'}{2}-1}.
\end{align*}
Consequently, we can derive as in (\ref{can}) that for any $t\leq 2$,
\begin{align*}
|\nabla^2_x\nabla_yT_tg(x,y)|&\leq C_2\Big(\|a\|_{C_b^{\delta,1}}+\|b\|_{C_b^{\delta,1}}\Big)\|g\|_{C_b^\delta}\int_0^{t}\!\!\int_{\mR^d}(t-s)^{-(d+2)/2}\\
&\qquad\qquad\times\exp\big(-c_0|z-x'|^2/(t-s)\big)\cdot|x'-z|^{\delta'}s^{\frac{\delta-\delta'}{2}-1}\dif x'\dif s\\
&\leq C_2\Big(\|a\|_{C_b^{\delta,1}}+\|b\|_{C_b^{\delta,1}}\Big)\|g\|_{C_b^\delta}\int_0^{t}(t-s)^{\delta'/2-1}s^{\frac{\delta-\delta'}{2}-1}\dif s\\
&\leq C_2t^{\delta/2-1}\Big(\|a\|_{C_b^{\delta,1}}+\|b\|_{C_b^{\delta,1}}\Big)\|g\|_{C_b^\delta},
\end{align*}
which yields (\ref{xy1}).
For $t>2$, we write
\begin{align*}
\nabla^2_x\nabla_yT_tg(x,y)&=\left(\int^t_{t-1}+\int^{t-1}_{1}+\int_0^{1}\right)\!\int_{\mR^d}\nabla^2_xp_{s}(x,z;y)\frac{\p\sL_0}{\p y}(z,y)T_{t-s}g(z,y)\dif z\dif s\\
&=:\cQ_1+\cQ_2+\cQ_3.
\end{align*}
Note that on $[t-1,t)$, we have $t-s\in(0,1]$. By (\ref{pp2}) and (\ref{can}), we can get
\begin{align*}
\cQ_1&\leq C_3\Big(\|a\|_{C_b^{\delta,1}}+\|b\|_{C_b^{\delta,1}}\Big)\|g\|_{C_b^\delta}\int_{t-1}^{t}\int_{\mR^d}\frac{1+|x|^m}{(1+s)^k(1+|z|^j)}(t-s)^{\delta/2-1}\dif z\dif s\\
&\leq C_3\Big(\|a\|_{C_b^{\delta,1}}+\|b\|_{C_b^{\delta,1}}\Big)\|g\|_{C_b^\delta}\frac{(1+|x|^m)}{(1+t)^{k}}.
\end{align*}
While for the second term, we can use (\ref{pp2}) and (\ref{can2}) to derive that
\begin{align*}
\cQ_2&\leq C_4\Big(\|a\|_{C_b^{\delta,1}}+\|b\|_{C_b^{\delta,1}}\Big)\|g\|_{C_b^\delta}\int_1^{t-1}\!\!\int_{\mR^d}\frac{1+|x|^m}{(1+s)^k(1+|z|^j)}\frac{1+|z|^m}{(1+(t-s))^k}\dif z\dif s\\
&\leq C_4\Big(\|a\|_{C_b^{\delta,1}}+\|b\|_{C_b^{\delta,1}}\Big)\|g\|_{C_b^\delta}\frac{(1+|x|^m)}{(1+t)^k}.
\end{align*}
To control the last term, we first claim that for every $k\in\mR_+$, there exist constants $C_5,m>0$ such that for any $x_1,x_2\in\mR^{d_1}$ and  $t\geq 1$,
\begin{align}\label{can3}
\sJ&:=\frac{\p\sL_0}{\p y}(x_1,y)T_{t}g(x_1,y)-\frac{\p\sL_0}{\p y}(x_2,y)T_{t}g(x_2,y)\no\\
&\leq C_5\Big(\|a\|_{C_b^{\delta,1}}+\|b\|_{C_b^{\delta,1}}\Big)\|g\|_{C_b^\delta}|x_1-x_2|^\delta\frac{1+|x_1|^m+|x_2|^m}{(1+t)^k}.
\end{align}
In fact, we can write
\begin{align*}
\sJ&\leq \bigg(\frac{\p\sL_0}{\p y}(x_1,y)-\frac{\p\sL_0}{\p y}(x_2,y)\bigg)T_{t}g(x_1,y)+\frac{\p\sL_0}{\p y}(x_2,y)\Big(T_{t}g(x_1,y)-T_{t}g(x_2,y)\Big)\\
&=:\sJ_1+\sJ_2.
\end{align*}
By the H\"older assumption on the coefficients and (\ref{can2}), it is easy to see that
$$
\sJ_1\leq C_6\Big(\|a\|_{C_b^{\delta,1}}+\|b\|_{C_b^{\delta,1}}\Big)\|g\|_{C_b^\delta}|x_1-x_2|^\delta\frac{1+|x_1|^m}{(1+t)^k}.
$$
On the other hand, we have by (\ref{ppp2}) that
\begin{align*}
\sJ_2&\leq C_7\Big(\|a\|_{C_b^{\delta,1}}+\|b\|_{C_b^{\delta,1}}\Big)\|g\|_{C_b^\delta}\sum_{\ell=1,2}\left|\int_{\mR^d}\big[\nabla^\ell_xp_t(x_1,x';y)-\nabla^\ell_xp_t(x_2,x';y)\big]g(x')\dif x'\right|\\
&=C_7\Big(\|a\|_{C_b^{\delta,1}}+\|b\|_{C_b^{\delta,1}}\Big)\|g\|_{C_b^\delta}\sum_{\ell=1,2}\bigg|\int_{\mR^d}\int_{\mR^d}\big[\nabla^\ell_xp_1(x_1,z;y)-\nabla^\ell_xp_1(x_2,z;y)\big]\\
&\qquad\qquad\qquad\qquad\qquad\qquad\qquad\quad\quad\times\big[p_{t-1}(z,x';y)-p_\infty(x';y)\big]g(x')\dif z\dif x'\bigg|\\
&\leq C_7\Big(\|a\|_{C_b^{\delta,1}}+\|b\|_{C_b^{\delta,1}}\Big)\|g\|_{C_b^\delta}|x_1-x_2|^\delta\!\int_{\mR^d}\int_{\mR^d}\!\Big(\exp(-c_0|x_1-z|^2)\\
&\qquad\qquad\qquad\qquad\qquad\qquad+\exp(-c_0|x_2-z|^2)\Big)\times\frac{1+|z|^m}{(1+t)^k(1+|x'|^j)}|g(x')|\dif z\dif x'\\
&\leq C_7\Big(\|a\|_{C_b^{\delta,1}}+\|b\|_{C_b^{\delta,1}}\Big)\|g\|_{C_b^\delta}|x_1-x_2|^\delta\frac{1+|x_1|^m+|x_2|^m}{(1+t)^k},
\end{align*}
where in the third inequality we also used (\ref{pin2}). Thus (\ref{can3}) is true.
It then follows by the same argument as in  (\ref{can}) that
\begin{align*}
\cQ_3
&\leq C_8\Big(\|a\|_{C_b^{\delta,1}}+\|b\|_{C_b^{\delta,1}}\Big)\|g\|_{C_b^\delta}\int^{1}_0\!\!\int_{\mR^d}s^{-(d+2)/2}\exp\big(-c_0|x-z|^2/2s\big)\\
&\qquad\qquad\qquad\qquad\qquad\qquad\qquad\quad\quad\times|x-z|^\delta\frac{1+|z|^m+|x|^m}{(1+(t-s))^k}\dif z\dif s\\
&\leq C_8\Big(\|a\|_{C_b^{\delta,1}}+\|b\|_{C_b^{\delta,1}}\Big)\|g\|_{C_b^\delta}\int^{1}_0\frac{1+|x|^m}{(1+(t-s))^k}s^{\delta/2-1}\dif s\\
&\leq C_8\Big(\|a\|_{C_b^{\delta,1}}+\|b\|_{C_b^{\delta,1}}\Big)\|g\|_{C_b^\delta}\frac{(1+|x|^m)}{(1+t)^k}.
\end{align*}
The proof is finished.
\end{proof}

We now establish the second order differentiability of $T_tg$ with respect to the $y$-variable. We have the following result.

\bl
Let {\bf (H$\sigma$)}, {\bf (H$b$)} and (\ref{cen}) hold. Assume that $a, b\in C_b^{\delta,2}$ and $g\in C_b^\delta$ with $0<\delta\leq 1$. Then we have
\begin{align}\label{d2}
\nabla^2_yT_tg(x,y)&=2\int_0^t\!\int_{\mR^d}p_{t-s}(x,z;y)\frac{\p\sL_0}{\p y}(z,y)\nabla_yT_{s}g(z,y)\dif z\dif s\no\\
&\quad+\int_0^t\!\!\int_{\mR^d}p_{t-s}(x,z;y)\frac{\p^2\sL_0}{\p y^2}(z,y)T_{s}g(z,y)\dif z\dif s.
\end{align}
Moreover, for any $0<t\leq 2$,
\begin{align}\label{es3}
|\nabla^2_yT_tg(x,y)|\leq C_0\Big[\big(\|a\|_{C_b^{\delta,1}}+\|b\|_{C_b^{\delta,1}}\big)^2+\big(\|a\|_{C_b^{\delta,2}}+\|b\|_{C_b^{\delta,2}}\big)\Big]\|g\|_{C_b^\delta},
\end{align}
and for any $k\in\mR_+$, there exists a constant $m>0$ such that for every $t>2$,
\begin{align}\label{es4}
|\nabla^2_yT_tg(x,y)|\leq C_0\Big[\big(\|a\|_{C_b^{\delta,1}}+\|b\|_{C_b^{\delta,1}}\big)^2+\big(\|a\|_{C_b^{\delta,2}}+\|b\|_{C_b^{\delta,2}}\big)\Big]\|g\|_{C_b^\delta}\frac{(1+|x|^m)}{(1+t)^k},
\end{align}
where $C_0>0$ is a constant depending only on $\lambda,d_1,d_2$ and $\|a\|_{C_b^{\delta,0}}, \|b\|_{C_b^{\delta,0}}$.
\el
\begin{proof}
The formula (\ref{d2}) has been proven in \cite[formula (34)]{P-V2}. Let us focus on estimates (\ref{es3}) and (\ref{es4}). In fact, for $t\leq2$, we can use (\ref{xy1}) and the same argument as in (\ref{es1}) to get that
\begin{align*}
|\nabla^2_yT_tg(x,y)|\leq C_0\Big(\|a\|_{C_b^{\delta,1}}+\|b\|_{C_b^{\delta,1}}\Big)^2\|g\|_{C_b^\delta}+\Big(\|a\|_{C_b^{\delta,2}}+\|b\|_{C_b^{\delta,2}}\Big)\|g\|_{C_b^\delta},
\end{align*}
which implies (\ref{es3}).
Now we prove the estimate (\ref{es4}).	To this end, we write
\begin{align*}
\nabla^2_yT_tg(x,y)&=2\bigg(\int_0^t\!\int_{\mR^d}p_{t-s}(x,z;y)\frac{\p\sL_0}{\p y}(z,y)\nabla_yT_{s}g(z,y)\dif z\dif s\\
&\quad-\int_0^\infty\!\int_{\mR^d}p_{\infty}(z;y)\frac{\p\sL_0}{\p y}(z,y)\nabla_yT_{s}g(z,y)\dif z\dif s\bigg)\\
&\quad+\bigg(\int_0^t\!\!\int_{\mR^d}p_{t-s}(x,z;y)\frac{\p^2\sL_0}{\p y^2}(z,y)T_{s}g(z,y)\dif z\dif s\\
&\quad-\int_0^\infty\!\int_{\mR^d}p_{\infty}(z;y)\frac{\p^2\sL_0}{\p y^2}(z,y)T_{s}g(z,y)\dif z\dif s\bigg)=:\cJ_1+\cJ_2.
\end{align*}
Note  by our assumption that $a, b\in C_b^{\delta,2}$, the second part $\cJ_2$ can be controlled in exactly the same way as in the estimate of $\nabla_yT_tg(x,y)$, i.e., we can get
$$
\cJ_2\leq C_1\Big(\|a\|_{C_b^{\delta,2}}+\|b\|_{C_b^{\delta,2}}\Big)\|g\|_{C_b^\delta}\frac{(1+|x|^m)}{(1+t)^k},
$$
where $C_1>0$ depends only on $\lambda,d_1,d_2$ and $\|a\|_{C_b^{\delta,0}}, \|b\|_{C_b^{\delta,0}}$.
Below we shall focus on the estimate of $\cJ_1$. As before, we write
\begin{align*}
\tfrac{1}{2}\cJ_1&= \int_0^{t/2}\!\!\int_{\mR^d}\big[p_{t-s}(x,z;y)-p_\infty(z;y)\big]\frac{\p\sL_0}{\p y}(z,y)\nabla_yT_sg(z,y)\dif z\dif s\\
&\quad+\int_{t/2}^t\int_{\mR^d}p_{t-s}(x,z;y)\frac{\p\sL_0}{\p y}(z,y)\nabla_yT_sg(z,y)\dif z\dif s\\
&\quad+\int_{t/2}^\infty\!\int_{\mR^d}p_\infty(z;y)\frac{\p\sL_0}{\p y}(z,y)\nabla_yT_sg(z,y)\dif z\dif s=:\cJ_{11}+\cJ_{12}+\cJ_{13}.
\end{align*}
For the first term, we have by (\ref{pin2}), (\ref{xy1})and (\ref{xy2}) that
\begin{align*}
\cJ_{11}&\leq C_1\Big(\|a\|_{C_b^{\delta,1}}+\|b\|_{C_b^{\delta,1}}\Big)^2\|g\|_{C_b^\delta}\int_0^{1}\!\!\int_{\mR^d}\frac{1+|x|^m}{(1+t-s)^k(1+|z|^j)}\dif z\dif s\\
&\quad+C_1\Big(\|a\|_{C_b^{\delta,1}}+\|b\|_{C_b^{\delta,1}}\Big)^2\|g\|_{C_b^\delta}\int_1^{t/2}\!\!\int_{\mR^d}\frac{1+|x|^m}{(1+t-s)^k(1+|z|^j)}\frac{1+|z|^m}{(1+s)^k}\dif z\dif s\\
&\leq C_1\Big(\|a\|_{C_b^{\delta,1}}+\|b\|_{C_b^{\delta,1}}\Big)^2\|g\|_{C_b^\delta}\frac{(1+|x|^m)}{(1+t)^k}.
\end{align*}
Using (\ref{pp1}) and (\ref{xy2}) again, we can control the second term by
\begin{align*}
\cJ_{12}
&\leq C\Big(\|a\|_{C_b^{\delta,1}}+\|b\|_{C_b^{\delta,1}}\Big)^2\|g\|_{C_b^\delta}\!\int^{1}_0\!\!\int_{\mR^d}\!s^{-d/2}\exp\big(\!-c_0|x-z|^2/2s\big)\frac{1+|z|^m}{(1+(t-s))^k}\dif z\dif s\\
&\quad+C\Big(\|a\|_{C_b^{\delta,1}}+\|b\|_{C_b^{\delta,1}}\Big)^2\|g\|_{C_b^\delta}\int^{t/2}_1\!\!\int_{\mR^d}\frac{1+|x|^m}{(1+|z|^{k})}\frac{1+|z|^m}{(1+(t-s))^k}\dif z\dif s\\
&\leq C\Big(\|a\|_{C_b^{\delta,1}}+\|b\|_{C_b^{\delta,1}}\Big)^2\|g\|_{C_b^\delta}\frac{(1+|x|^m)}{(1+t)^k}.
\end{align*}
Finally, we have by (\ref{pin1}) and (\ref{xy2}) that
	\begin{align*}
	\cI_3&\leq C\Big(\|a\|_{C_b^{\delta,1}}+\|b\|_{C_b^{\delta,1}}\Big)^2\|g\|_{C_b^\delta}\int_{t/2}^\infty\!\int_{\mR^d}\frac{1}{1+|z|^j}\cdot\frac{1+|z|^m}{(1+s)^k}\dif z\dif s\\
	&\leq C\Big(\|a\|_{C_b^{\delta,1}}+\|b\|_{C_b^{\delta,1}}\Big)^2\|g\|_{C_b^\delta}\frac{(1+|x|^m)}{(1+t)^k}.
	\end{align*}
The proof is finished.
\end{proof}

With the above preparations, we can establish the following regularity of $T_tf$ with respect to the parameter $y$.

\bl\label{dyf}
Let {\bf (H$\sigma$)}, {\bf (H$b$)} and (\ref{cen}) hold. Assume $a, b\in C_b^{\delta,\ell}$ and $f\in C_b^{\delta,\ell}$ with $0<\delta\leq 1$ and $\ell=1,2$. Then we have:\\
(i) (Case $\ell=1$ and $0<t\leq 2$):
\begin{align*}
|\nabla_y T_tf(x,y)|&\leq C_0\Big[\big(\|a\|_{C_b^{\delta,0}}+\|b\|_{C_b^{\delta,0}}\big)\|f\|_{C_b^{\delta,1}}+\big(\|a\|_{C_b^{\delta,1}}+\|b\|_{C_b^{\delta,1}}\big)\|f\|_{C_b^{\delta,0}}\Big];
\end{align*}
(ii) (Case $\ell=1$ and $t>2$): for any $k\in\mR_+$, there exists a constant $m>0$ such that for every $t>2$,
\begin{align*}
|\nabla_y T_tf(x,y)|\leq C_0&\Big[\big(\|a\|_{C_b^{\delta,0}}+\|b\|_{C_b^{\delta,0}}\big)\|f\|_{C_b^{\delta,1}}\no\\
&+\big(\|a\|_{C_b^{\delta,1}}+\|b\|_{C_b^{\delta,1}}\big)\|f\|_{C_b^{\delta,0}}\Big]\frac{(1+|x|^m)}{(1+t)^k};
\end{align*}
(iii) (Case $\ell=2$ and $0<t\leq 2$):
\begin{align}\label{es6}
|\nabla_y^2 T_tf(x,y)|\leq &C_0\Big[\big(\|a\|_{C_b^{\delta,0}}+\|b\|_{C_b^{\delta,0}}\big)\|f\|_{C_b^{\delta,2}}+\big(\|a\|_{C_b^{\delta,1}}+\|b\|_{C_b^{\delta,1}}\big)\|f\|_{C_b^{\delta,1}}\no\\
&+\Big(\big(\|a\|_{C_b^{\delta,1}}+\|b\|_{C_b^{\delta,1}}\big)^2+\big(\|a\|_{C_b^{\delta,2}}+\|b\|_{C_b^{\delta,2}}\big)\Big)\|f\|_{C_b^{\delta,0}}\Big];
\end{align}
(iv) (Case $\ell=2$ and $t>2$): for any $k\in\mR_+$, there exists a constant $m>0$ such that  for every $t>2$,
\begin{align}\label{es66}
|\nabla_y^2& T_tf(x,y)|\leq C_0\Big[\big(\|a\|_{C_b^{\delta,0}}+\|b\|_{C_b^{\delta,0}}\big)\|f\|_{C_b^{\delta,2}}+\big(\|a\|_{C_b^{\delta,1}}+\|b\|_{C_b^{\delta,1}}\big)\|f\|_{C_b^{\delta,1}}\no\\
&+\Big(\big(\|a\|_{C_b^{\delta,1}}+\|b\|_{C_b^{\delta,1}}\big)^2+\big(\|a\|_{C_b^{\delta,2}}+\|b\|_{C_b^{\delta,2}}\big)\Big)\|f\|_{C_b^{\delta,0}}\Big]\frac{(1+|x|^m)}{(1+t)^k},
\end{align}
where $C_0>0$ is a constant depending only on $\lambda,d_1,d_2$ and $\|a\|_{C_b^{\delta,0}}, \|b\|_{C_b^{\delta,0}}$.
\el

\begin{proof}
We only prove the above estimates when $\ell=2$. The corresponding estimates	for $\ell=1$ follows by the same arguments. In fact, we have
	\begin{align*}
	\nabla^2_yT_tf(x,y)&=\sum_{\ell=0}^2 C^\ell_2\nabla_y^\ell T_tg(x,y)\Big|_{g=\nabla_y^{2-\ell}f}\\
	&=T_tg(x,y)\Big|_{g=\nabla^2_yf}+2\nabla_yT_t\nabla_yg(x,y)\Big|_{g=\nabla_yf}+\nabla^2_yT_tg(x,y)\Big|_{g=f}\\
	&=:\cK_1+\cK_2+\cK_3.
	\end{align*}
	When $t\leq 2$, it is obvious that
	$$
	\cK_1\leq C_1\Big(\|a\|_{C_b^{\delta,0}}+\|b\|_{C_b^{\delta,0}}\Big)\|\nabla_y^2f\|_{C_b^{\delta,0}}\leq C_1\Big(\|a\|_{C_b^{\delta,0}}+\|b\|_{C_b^{\delta,0}}\Big)\|f\|_{C_b^{\delta,2}}.
	$$
	For the second term, we have by (\ref{es1}) that
	$$
	\cK_2\leq C_2\Big(\|a\|_{C_b^{\delta,1}}+\|b\|_{C_b^{\delta,1}}\Big)\|\nabla_yf\|_{C_b^{\delta,0}}\leq C_2\Big(\|a\|_{C_b^{\delta,1}}+\|b\|_{C_b^{\delta,1}}\Big)\|f\|_{C_b^{\delta,1}}.
	$$
	Finally, using (\ref{es3}) we can control the third term by
	\begin{align*}
	\cK_3\leq C_3\Big[\big(\|a\|_{C_b^{\delta,1}}+\|b\|_{C_b^{\delta,1}}\big)^2+\big(\|a\|_{C_b^{\delta,2}}+\|b\|_{C_b^{\delta,2}}\big)\Big]\|f\|_{C_b^{\delta,0}},
	\end{align*}
	which in turn yields (\ref{es6}).
The estimate (\ref{es66}) can be proved similarly by replacing (\ref{es1}) and (\ref{es3}) with (\ref{es2}) and (\ref{es4}), respectively.  The proof is finished.
\end{proof}

\section{Strong convergence with order $(\alpha\wedge1)/2$}

In  this section, we study the strong convergence of the multi-scale system (\ref{sde0}) to the effective equation (\ref{sde1}). To this end, we assume that
$$
G(t,x,y)\equiv G(t,y),
$$
i.e., the diffusion coefficient $G$ in the slow equation does not dependent on the $x$-variable. Note that in this case, we have
$$
\bar G(t,y)=G(t,y).
$$
We shall always assume {\bf (H$\sigma$)}, {\bf (H$_G$)}, {\bf (H$b$)} hold, and that the coefficients $a$ and $b$ are H\"older continuous with respect to $x$ uniformly in $y$, and that the coefficient $G$ is H\"older continuous with respect to $y$ uniformly in $t$.

\subsection{Zvonkin transform}
Due to the low regularity assumptions on the coefficients of the system (\ref{sde0}), it is not possible to prove the strong convergence of $Y_t^\eps$ to $\bar Y_t$ directly. For this reason, we shall use Zvonkin's argument to transform the equations for $Y_t^\eps$ and $\bar Y_t$ into new ones.
Let us first prove the following regularity result for the averaged drift coefficient.

\bl\label{hff}
Assume that $a, b\in C_b^{\delta,\alpha}$ and $F\in C_b^{\alpha/2,\delta,\alpha}$ with $0<\delta,\alpha\leq 1$. Let $\bar F$ be defined as in (\ref{bf}). Then we have $\bar F\in C_b^{\alpha/2,\alpha}$.
\el
\begin{proof}
The $\alpha/2$-H\"older continuity with respect to the $t$ variable follows directly by the definition of $\bar F$. Let us prove the H\"older continuous with respect to $y$.	We write for $y_1,y_2\in\mR^{d_2}$
\begin{align*}
\bar F(t,y_1)-\bar F(t,y_2)&=\int_{\mR^{d_1}}\![F(t,x,y_1)-F(t,x,y_2)]\mu^{y_1}(\dif x)\\
&\quad+\int_{\mR^{d_1}}F(t,x,y_2)\big[\mu^{y_1}(\dif x)-\mu^{y_2}(\dif x)\big]=:\sK_1+\sK_2.
\end{align*}
It is easy to see that there exists a constant $C_1>0$ such that
$$
\sK_1\leq C_1\big(|y_1-y_2|^\alpha\wedge1\big).
$$
For the second term, by the same argument as in (\ref{d1}), we  get
\begin{align*}
	\sK_2&=\int_{\mR^d}F(t,x',y)\big[p_\infty(x';y_1)\dif x'-p_\infty(x';y_2)\dif x'\big]\\
	&=\lim_{t\to\infty}\int_0^t\!\!\int_{\mR^d}p_{t-s}(x,z;y)[\sL_0(z,y_1)-\sL_0(z,y_2)]\left(\int_{\mR^d}p_s(z,x';y)F(t,x',y)\dif x'\right)\dif z\dif s\\
	&=\int_0^\infty\!\!\int_{\mR^d}p_\infty(z;y)[\sL_0(z,y_1)-\sL_0(z,y_2)]\left(\int_{\mR^d}p_s(z,x';y)F(t,x',y)\dif x'\right)\dif z\dif s.
	\end{align*}
Thus, we have by (\ref{pin1}), (\ref{can}) and (\ref{can2}) that
	\begin{align*}
	\sK_2
	&\leq C_2\big(|y_1-y_2|^\alpha\wedge1\big)\bigg(\int_0^2\!\!\int_{\mR^d}\frac{1}{(1+|z|^j)}s^{\delta/2-1}\dif z\dif s\\
	&\quad+\int_2^\infty\!\!\!\int_{\mR^d}\frac{1}{(1+|z|^j)}\frac{(1+|z|)^m}{(1+s)^k}\dif z\dif s\bigg)\leq C_2\big(|y_1-y_2|^\alpha\wedge1\big),
	\end{align*}
where $C_2>0$ is a constant.
	The proof is finished.
\end{proof}

Below, we shall fix a $T>0$ to be sufficiently small. Recall that $\bar\sL$ is defined by (\ref{lby}).
Consider the following backward PDE in $\mR^{d_2}$:
\begin{equation}\left\{\begin{array}{l}\label{PDE2}
\displaystyle
\p_tv(t,y)+\bar{\mathscr{L}}v(t,y)+\bar F(t,y)=0,\quad t\in [0, T),\\
v(T,y)=0.
\end{array}\right.
\end{equation}
Under our assumptions on the coefficients and by Lemma \ref{hff}, it is well known that there exits a unique solution $v\in L^\infty\big([0,T];C^{2+\alpha}_{b}(\mR^{d_2})\big)\cap C_b^{1+\alpha/2}\big([0,T];L^\infty(\mR^{d_2})\big)$ for equation \eqref{PDE2}. Moreover, we can choose $T$ small enough so that for any $0<t<T$,
$$
1/2\leq |\nabla_yv(t,y)|\leq 2,\quad\forall y\in\mR^{d_2}.
$$
Define the transformed function by
$$
\Phi(t,y):=y+v(t,y).
$$
Then, the map $y\rightarrow\Phi(t,y)$ forms a $C^1$-diffeomorphism and
\begin{align}\label{dd}
1/2\leq\|\nabla_y\Phi\|_\infty\leq 3.
\end{align}

Now, let us define the new processes by
\begin{align}\label{vv}
\bar V_t:=\Phi(t,\bar{Y}_t)\quad\text{and}\quad V^{\eps}_t:=\Phi(t,Y^{\eps}_t).
\end{align}
We have the following result.
\bl
Let $\bar V_t$ and $V_t^\eps$ be defined by (\ref{vv}). Then we have
\begin{align}\label{z1}
\dif \bar V_t=G(t,\bar Y_t)\nabla_y \Phi(t,\bar{Y}_t)\dif W^2_t,\quad V_0=\Phi(0, y)
\end{align}
and
\begin{align}\label{z2}
\dif V^{\eps}_t=\big[F(t,X^{\eps}_t, Y^{\eps}_t)&-\bar{F}(t,Y^{\eps}_t)\big]\nabla_y \Phi(t,Y^{\eps}_t)\dif t\no\\
&+G(t,Y_t^\eps)\nabla_y \Phi(t,Y^{\eps}_t)\dif W^2_t,\quad V_0^\eps=\Phi(0, y).
\end{align}
\el
\begin{proof}
	Using It\^o's formula, we have
	\begin{align*}
	v(t,Y_t^\eps)&=v(0,y)+\int_0^t\big(\p_s+\sL_1\big)v(s,Y_s^\eps)\dif s+\int_0^tG(s,Y_s^\eps)\nabla_yv(s,Y_s^\eps)\dif W_s^2\\
	&=v(0,y)+\int_0^t\big(\p_s+\bar\sL\big)v(s,Y_s^\eps)\dif s+\int_0^tG(s,Y_s^\eps)\nabla_yv(s,Y_s^\eps)\dif W_s^2\\
	&\quad+\int_0^t\big[F(s,X^{\eps}_s, Y^{\eps}_s)-\bar{F}(s,Y^{\eps}_s)\big]\nabla_y v(s,Y^{\eps}_s)\dif s\\
	&=v(0,y)-\int_0^t\bar F(s,Y_s^\eps)\dif s+\int_0^tG(s,Y_s^\eps)\nabla_yv(s,Y_s^\eps)\dif W_s^2\\
	&\quad+\int_0^t\big[F(s,X^{\eps}_s, Y^{\eps}_s)-\bar{F}(s,Y^{\eps}_s)\big]\nabla_y v(s,Y^{\eps}_s)\dif s,
	\end{align*}
	where in the last equality we used (\ref{PDE2}). This together with the equation for $Y_t^\eps$ yields (\ref{z2}). The proof of (\ref{z1}) is easier and follows by the same argument.
\end{proof}

\subsection{Proof of Theorem \ref{main1}}

We  first prepare the following mollifying approximation result.
For simplification, let us set
\begin{align}\label{hf}
\hat F(t,x,y):=\big[F(t,x,y)-\bar{F}(t,y)\big]\nabla_y\Phi(t,y).
\end{align}
Let $\rho_1:\mR\to[0,1]$ and $\rho_2:\mR^{d_2}\to[0,1]$ be two smooth radial convolution kernel functions
such that
$\int_\mR\rho_1(r)\dif r=\int_{\mR^{d_2}}\rho_2(y)\dif y=1$, and for any $k\geq 1$, $|\nabla^k\rho_1|\leq C_k\rho_1(x)$ and  $|\nabla^k\rho_2|\leq C_k\rho_2(x)$, where $C_k>0$ are constants.
%%% ($\rho(x)=C_0\e^{-\sqrt{1+|x|^2}}, for instance)$
For every $n\in\mN^*$, set
$$
\rho_1^n(r):=n^2\rho_1(n^2r)\quad \text{and}\quad
\rho_2^n(y):=n^{d_2}\rho_2(ny).$$
We define the mollifying approximations of $\hat F$  by
\begin{align}\label{fn}
\hat F_n(t,x,y):=\int_{\mR^{d_2+1}}\hat F(t-s,x,y-z)\rho_2^{n}(z)\rho_1^n(s)\dif z\dif s.
\end{align}
Similarly, we define the mollifying approximations of $a, b$  by
\begin{align}\label{bn}
a_n(x,y):=\int_{\mR^{d_2}}a(x,y-z)\rho_2^{n}(z)\dif z,\quad b_n(x,y):=\int_{\mR^{d_2}}b(x,y-z)\rho_2^{n}(z)\dif z.
\end{align}
We have the following easy result, which will play important role below.

\bl\label{momo}
Assume that $a,b\in C_b^{\delta,\alpha}$ and $F\in C_b^{\alpha/2,\delta,\alpha}$ with $0<\delta,\alpha\leq1$. Then we have
\begin{align}\label{n1}
\|\hat F-\hat F_n\|_\infty+\|a-a_n\|_\infty+\|b-b_n\|_\infty\leq C_0n^{-\alpha},
\end{align}
and
\begin{align}\label{n2}
\|\hat F_n\|_{C_b^{1,\delta,\alpha}}+\|\hat F_n\|_{C_b^{\alpha,\delta,2}}+ \|a_n\|_{C_b^{\delta,2}}+\|b_n\|_{C_b^{\delta,2}}\leq C_0n^{2-\alpha},
\end{align}
where $C_0>0$ is a constant independent of $n$.
\el
\begin{proof}
According to Lemma \ref{hff}, it is easy to check that $\hat F\in C_b^{\alpha/2,\delta,\alpha}$.  By the definition of $\hat F_n$, we have
	\begin{align*}
	|\hat F(t,x,y)-\hat F_n(t,x,y)|&\leq \int_{\mR^{d_2+1}}\big|\hat F(t-s,x,y-z)-\hat F(t,x,y)\big|\cdot\rho_2^{n}(z)\rho_1^n(s)\dif z\dif s\\
	&\leq C_1\int_{\mR^{d_2+1}}\big(s^{\alpha/2}+|z|^{\alpha}\big)\cdot\rho_2^{n}(z)\rho_1^n(s)\dif z\dif s\leq C_1n^{-\alpha}.
	\end{align*}
Furthermore, we have
\begin{align*}
	|\p_t\hat F_n(t,x,y)|&\leq \int_{\mR^{d_2+1}}\big|\hat F(t-s,x,y-z)-\hat F(t,x,y-z)\big|\cdot|\rho_2^{n}(z)|\p_s\rho_1^n(s)\dif z\dif s\\
	&\leq C_2n^2\int_{\mR^{d_2+1}}s^{\alpha/2}\rho_2^{n}(z)\rho_1^n(s)\dif z\dif s\leq C_2n^{2-\alpha},
	\end{align*}
	and
	\begin{align*}
	|\nabla_y^2\hat F_n(t,x,y)|&\leq \int_{\mR^{d_2+1}}\big|\hat F(t-s,x,y-z)-\hat F(t-s,x,y)\big|\cdot|\nabla_z^2\rho_2^{n}(z)|\rho_1^n(s)\dif z\dif s\\
	&\leq C_2n^2\int_{\mR^{d_2+1}}|z|^{\alpha}\rho_2^{n}(z)\rho_1^n(s)\dif z\dif s\leq C_2n^{2-\alpha}.
	\end{align*}
The other estimates can be proved similarly.
\end{proof}

Now, we are in the position to give:

\begin{proof}[Proof of Theorem \ref{main1}]
	Let us first assume that $T>0$ is sufficiently small so that (\ref{dd}) holds.
	As a result, we have for any $t\in [0, T]$,
	\begin{align}\label{aaa}
	\mE\big|Y_{t}^\eps-\bar Y_{t}\big|^2\leq 2\mE\big|V_{t}^\eps-\bar V_{t}\big|^2.
	\end{align}
	Hence, we shall focus on the convergence of $V_t^\eps$ to $\bar V_t$. Recall the definition of $\hat F$ by (\ref{hf}), and let $\hat F_n$ be given by (\ref{fn}). According to (\ref{z1}) and (\ref{z2}), we write
	\begin{align*}
	V_t^\eps-\bar V_t&=\int^t_0\big[G(s,Y_s^\eps)\nabla_y \Phi(s,Y^{\eps}_s)-G(s,\bar Y_s)\nabla_y \Phi(s,\bar Y_s)\big]\dif W^2_s\\
	&\quad+\int^t_0\big[\hat F(s,X^{\eps}_s, Y^{\eps}_s)-\hat F_n(s,X^{\eps}_s, Y^{\eps}_s)\big]\dif s+\int_0^{t}\hat F_n(s,X^{\eps}_s, Y^{\eps}_s)\dif s.
	\end{align*}
	Thus, taking expectation and using  Burkholder-Davis-Gundy's inequality we can get that there exists a $C_0>0$ such that
	\begin{align*}
	\mE|V_{t}^\eps-V_{t}|^2&\leq C_0\mE\left(\int^{t}_0\big|G(s,Y_s^\eps)\nabla_y \Phi(s,Y^{\eps}_s)-G(s,\bar Y_s)\nabla_y \Phi(s,\bar Y_s)\big|^2\dif s\right)\\
	&\quad+C_0\mE\left|\int_0^{t}\!\!\big[\hat F(s,X^{\eps}_s, Y^{\eps}_s)-\hat F_n(s,X^{\eps}_s, Y^{\eps}_s)\big]\dif s\right|^2\\
	&\quad+C_0\mE\left|\int_0^{t}\hat F_n(s,X^{\eps}_s, Y^{\eps}_s)\dif s\right|^2=:\sQ_1(t,\eps)+\sQ_2(t,\eps)+\sQ_3(t,\eps).
	\end{align*}
	Below, we divide the proof into three steps to control each term on the right hand side separately.
	
	\vspace{1mm}
	\noindent
	{\bf Step 1} (Control of $\sQ_1(t,\eps)$).
	Note that the function
	$$
	y\to G(t,\cdot)\nabla_y\Phi(t,\cdot)\in  C_b^1(\mR^{d_2}).
	$$
	As a result, we easily have that
	\begin{align}\label{ae}
	\sQ_1(t,\eps)\leq C_1\mE\left(\int^{t}_0|Y_s^\eps-\bar Y_s|^2\dif s\right),
	\end{align}
	where $C_1>0$ is a constant independent of $\eps$.
	
	\vspace{1mm}
	\noindent
	{\bf Step 2} (Control of $\sQ_2(t,\eps)$).
	The estimate of this term follows by an easy consequence of (\ref{n1}), which in turn yields that
	\begin{align}\label{p1}
	\sQ_{2}(t,\eps)\leq C_2\|\hat F-\hat F_n\|^2_\infty\leq C_2n^{-2\alpha},
	\end{align}
	where $C_2$ is a positive constant independent of $n$ and $\eps$.
	
\vspace{1mm}
	\noindent
	{\bf Step 3} (Control of $\sQ_3(t,\eps)$).
	We use the technique of the Poisson equation to control the third part. Let $a_n, b_n$ be defined by (\ref{bn}), and denote by $\sL_0^n(x,y)$ the operator $\sL_0(x,y)$ with coefficients $a, b$ replaced by $a_n, b_n$, i.e.,
	\begin{align}\label{l0n}
	\sL_0^n(x,y):=\sum_{i,j} a^{ij}_n(x,y)\frac{\p^2}{\p x_i\p x_j}+b_n(x,y)\cdot\nabla_x.
	\end{align}
	Let $\Psi_n$ be the solution to the following Poisson equation in $\mR^{d_1}$:
	$$
	\sL_0^n(x,y)\Psi_n(t,x,y)=\hat F_n(t,x,y),
	$$
	where $(t,y)\in\mR_+\times\mR^{d_2}$ are viewed as parameters.
Note that $\hat F_n$ satisfies the centering condition (\ref{cen}). Thus, according to Theorem \ref{popde},
	we can use  It\^o's formula to get that for any $t>0$,
	\begin{align*}
	&\Psi_n({t},X_{t}^\eps,Y_{t}^\eps)=\Psi_n(0,x,y)+\int_0^{t}\big(\p_s+\sL_1\big)\Psi_n(s,X_s^\eps,Y_s^\eps)\dif s\\
	&\qquad\qquad\qquad\qquad\qquad+\int_0^{t}\frac{1}{\eps}\sL_0\Psi_n(s,X_s^\eps,Y_s^\eps)\dif s+\frac{1}{\sqrt{\eps}}M^1_{t}+M^2_{t},
	\end{align*}
	where $\sL_1$ is given by (\ref{l1}), and  for $i=1,2$, $M^i_{t}$ are martingales defined by
	$$
	M^1_{t}:=\int_0^t\nabla_x\Psi_n(s,X_s^\eps,Y_s^\eps)\sigma(s,X_s^\eps)\dif W_s^1,
	$$
	and
	$$
	M^2_{t}:=\int_0^t\nabla_y\Psi_n(s,X_s^\eps,Y_s^\eps)G(s,Y_s^\eps)\dif W_s^2.
	$$
	This in turn yields that
	\begin{align*}
	\int^{t}_0\hat F_n(s,X^{\eps}_s, Y^{\eps}_s)\dif s
	&=\eps\Psi_n({t},X_{t}^\eps,Y_{t}^\eps)-\eps\Psi_n(0,x,y)-\sqrt{\eps}M^1_{t}-\eps M^2_{t}\\
	&\quad+\int_0^{t}\big[b_n(X_s^\eps,Y_s^\eps)-b(X_s^\eps,Y_s^\eps)\big]\nabla_x\Psi_n(s,X_s^\eps,Y_s^\eps)\dif s\\
&\quad+\int_0^{t}\big[a_n(X_s^\eps,Y_s^\eps)-a(X_s^\eps,Y_s^\eps)\big]\nabla^2_x\Psi_n(s,X_s^\eps,Y_s^\eps)\dif s\\
	&\quad-\eps\int_0^{t}\big(\p_s+\sL_1\big)\Psi_n(s,X_s^\eps,Y_s^\eps)\dif s.
	\end{align*}
	Taking this back into the definition of $\sQ_{3}(t,\eps)$ and by (\ref{key1}), we have that there exists a constant $m>0$ such that
	\begin{align*}
	&\sQ_{3}(t,\eps)
	\leq C_3\bigg[\eps^2\mE(1+|X^{\eps}_t|^{2m})+\eps\mE|M^1_{t}|^2+\eps^2\mE|M^2_{t}|^2\\ &+\mE\left(\int_0^{t}\!\Big(\big|b_n(X_s^\eps,Y_s^\eps)-b(X_s^\eps,Y_s^\eps)\big|^2\!+\!\big|a_n(X_s^\eps,Y_s^\eps)-a(X_s^\eps,Y_s^\eps)\big|^2\Big) (1+|X^{\eps}_s|^{2m})\dif s\!\right)\\
	&+\eps^2\mE\left|\int_0^{t}\big(\p_s+\sL_y\big)\Psi_n(s,X_s^\eps,Y_s^\eps)\dif s\right|^2\bigg]=:\sQ_{31}(t,\eps)+\sQ_{32}(t,\eps)+\sQ_{33}(t,\eps).
	\end{align*}
Note that the assumptions {\bf (H$\sigma$) } and {\bf (H$b$)} hold uniformly in $y$. Hence, it follows by \cite[Lemma 1]{Ve5} (see also \cite[Lemma 2]{P-V2}) that for any $k>0$,
\begin{align}
\mE|X^{\eps}_t|^{k}\leq C(1+|x|^{k}),\label{ME}
\end{align}
where $C$ is a positive constant independent of $\eps$. As a result, we can control the first term by (\ref{key3}) and (\ref{n2}) that
\begin{align*}
\sQ_{31}(t,\eps)&\leq C_4\Big(\eps+\eps^2\mE\int_0^t\|\nabla_y\Psi(\cdot,X_s^\eps,\cdot)\|_\infty^2\dif s\Big)\\
&\leq C_4\Big(\eps+\eps^2\Big(\|a_n\|_{C_b^{\delta,1}}+\|b_n\|_{C_b^{\delta,1}}+\|\hat F_n\|_{C_b^{0,\delta,1}}\Big)\mE\left(\int_0^{t}(1+|X^{\eps}_s|^{2m})\dif s\right)\\
&\leq C_4 \Big(\eps+\eps^2n^{2(1-\alpha)}\Big).
\end{align*}
For the second term, by (\ref{n1}) and (\ref{ME}), it is easy to see that
	\begin{align*}
	\sQ_{32}(t,\eps)&\leq C_5\Big(\|b_n-b\|_\infty^2+\|a_n-a\|_\infty^2\Big)\mE\left(\int_0^{t}(1+|X^{\eps}_s|^{2m})\dif s\right)\leq  C_5n^{-2\alpha}.
	\end{align*}
To estimate the last part, we first note that by (\ref{pro}) and viewing $t$ as a parameter,we have that for any $s>0$, $x\in\mR^{d_1}$ and $y\in\mR^{d_2}$,
\begin{align*}
|\p_s\Psi_n(s,x,y)|\leq C_6\|\p_s\hat F_n\|_{C_b^{0,\delta,0}}(1+|x|^m)\leq C_6n^{2-\alpha}(1+|x|^m),
\end{align*}
where the last inequality follows by (\ref{n2}).
On the other hand,  by reviewing $y$ as a parameter, we have by (\ref{key2}) and (\ref{n2}) that
\begin{align*}
	\|\sL_y\Psi_n(\cdot,x,\cdot)\|_\infty\!&\leq \!C_6\Big[\big(\|a_n\|_{C_b^{\delta,0}}+\|b_n\|_{C_b^{\delta,0}}\big)\|\hat F_n\|_{C_b^{0,\delta,2}}\!+\!\big(\|a_n\|_{C_b^{\delta,1}}+\|b_n\|_{C_b^{\delta,1}}\big)\|\hat F_n\|_{C_b^{0,\delta,1}}\no\\
&+\Big(\big(\|a_n\|_{C_b^{\delta,1}}+\|b_n\|_{C_b^{\delta,1}}\big)^2+\big(\|a_n\|_{C_b^{\delta,2}}+\|b_n\|_{C_b^{\delta,2}}\big)\Big)\|\hat F_n\|_{C_b^{0,\delta,0}}\Big](1+|x|^m)\\
&\leq C_6\Big(n^{2-\alpha}+n^{2-2\alpha}\Big)(1+|x|^m)\leq C_6n^{2-\alpha}(1+|x|^m).
\end{align*}
	As a result, we have
	$$
	\sQ_{33}(t,\eps)\leq C_7\eps^2n^{2(2-\alpha)}.
	$$
	Combing the above estimates, we get
	$$
	\sQ_3(t,\eps)\leq C_8\big(\eps+n^{-2\alpha}+\eps^2n^{2(2-\alpha)}\big).
	$$
	
	Now, in view of (\ref{aaa}), (\ref{ae}) and (\ref{p1}), we arrive at
	$$
	\mE\big|Y_{t}^\eps-\bar Y_{t}\big|^2\leq C_9\mE\left(\int^{t}_0|Y_s^\eps-\bar Y_s|^2\dif s\right)+C_9\Big(n^{-2\alpha}+\eps+\eps^2n^{2(2-\alpha)}\Big).
	$$
	Taking $n=\eps^{-1/2}$, we get
	$$
	\mE\big|Y_{t}^\eps-\bar Y_{t}\big|^2\leq C_9\mE\left(\int^{t}_0|Y_s^\eps-\bar Y_s|^2\dif s\right)+C_9\eps^{\alpha\wedge1},
	$$
	which in turn yields by Gronwall's inequality that
	$$
	\mE\big|Y_{t}^\eps-\bar Y_{t}\big|^2\leq C_T\eps^{\alpha\wedge1}.
	$$
	For general $T>0$, the result can be proved by induction and analogous arguments.
	So, the whole proof is finished.
\end{proof}

\subsection{Proof of Theorem \ref{main4}}

We use Theorem \ref{main1} and the radial  truncation technique to prove Theorem \ref{main4}.

\begin{proof}[Proof of Theorem \ref{main4}] For each $n\in\mN$, define
the new coefficients by
\begin{eqnarray*}
b_n(x,y):=\left\{ \begin{aligned}
&b(x,y),\qquad\quad\!\! |y|\leq n,\\
&b(x,ny/|y|)\quad |y|> n,
\end{aligned} \right.,\quad \sigma_n(x,y):=\left\{ \begin{aligned}
&\sigma(x,y),\qquad\,\,\,\, |y|\leq n,\\
&\sigma(x,ny/|y|)\quad |y|> n,
\end{aligned} \right.
\end{eqnarray*}
and
\begin{eqnarray*}
F_n(t,x,y):=\left\{ \begin{aligned}
&F(t,x,y),\qquad\,\,\,\, |y|\leq n,\\
&F(t,x,ny/|y|)\quad |y|> n,
\end{aligned} \right.,\quad G_n(t,y):=\left\{ \begin{aligned}
&G(t,y),\qquad\,\,\,\, |y|\leq n,\\
&G(t,ny/|y|)\quad |y|> n.
\end{aligned} \right.
\end{eqnarray*}
It is easy to check that $b_n, \sigma_n, F_n, G_n$ satisfy the conditions in Theorem \ref{main1}.
Let $(X^{n,\eps}_t,Y^{n,\eps}_t)$ be the solution to SDE (\ref{sde0}) with coefficients $b, \sigma, F, G$ replaced by $b_n, \sigma_n, F_n, G_n$. Then for any $T>0$, we have by Theorem \ref{main1} that
$$
\sup_{t\in[0, T]}\mE\big|Y_{t}^{n,\eps}-\bar Y_{t}^n\big|^2\to0\quad\text{as}\quad\eps\to0,
$$
where $\bar Y_t^n$ is the solution of the following new averaged equation:
\begin{align*}
\dif \bar Y_t^n=\bar F_n(t,\bar Y^n_t)\dif t+G_n(t,\bar Y^n_t)\dif W^2_t,\quad\bar Y^n_0=y.
\end{align*}
Here, $\bar F_n(t,y):=\int_{\mR^{d_1}} F_n(t,x,y)\mu^{y}_n(dx)$, and $\mu^{y}_n(\dif x)$ is the unique invariant measure of the transition semigroup of the following frozen equation:
\begin{align*}
\dif X_t^{n,y}=b_n(X_t^{n,y},y)\dif t+\sigma_n(X_t^{n,y},y)\dif W_t^1,\quad X_0^{n,y}=x.
\end{align*}
For every $\eps>0$, define the stopping time by
$$
\tau_n^\eps:=\inf\{t\geq 0: |Y^{\eps}_t|+|\bar{Y}_t|\geq n \}.
$$
Then, by the construction of the new coefficients and  the uniqueness of the strong solution to SDE (\ref{sde0}), it holds
$$
Y^\eps_t=Y^{n,\eps}_t,\quad \forall t\in[0,\tau_{n}^\eps].
$$
On the other hand, note that for every $|y|\leq n$, we also have $\mu^{y}_n(dx)=\mu^{y}(dx)$. This implies that for $|y|\leq n$,
\begin{align*}
\bar F_n(t, y)=\int_{\mR^{d_1}} F_n(t,x,y)\mu^{y}_n(dx)=\int_{\mR^{d_1}} F(t,x,y)\mu^{y}(dx)=\bar F(t, y),
\end{align*}
which together with the uniqueness of the strong solution to  SDE (\ref{sde1}) means
$$
\bar Y^n_t=\bar Y_t,\quad \forall t\in[0,\tau_{n}^\eps].
$$
As a result, we can deduce that for some $\beta>2$ and $C>0$,
\begin{align*}
\sup_{t\in[0, T]}&\mE\big|Y_{t}^\eps-\bar Y_{t}\big|^2\leq \sup_{t\in[0, T]}\mE\big(|Y_{t}^\eps-\bar Y_{t}|^2\cdot1_{\{t\leq \tau^\eps_n\}}\big)+\sup_{t\in[0, T]}\mE\big(|Y_{t}^\eps-\bar Y_{t}|^2\cdot1_{\{t>\tau^\eps_n\}}\big)\\
&\leq  \sup_{t\in[0, T]}\mE|Y^{n,\eps}_{t}-\bar Y^n_{t}|^2+C\left[\sup_{t\in[0, T]}\mE\big(|Y_t^\eps|^\beta+|\bar Y_t|^\beta\big)\right]^{2/\beta}\left[\mP(T>\tau^\eps_n)\right]^{(\beta-2)/\beta}\\
&\leq \sup_{t\in[0, T]}\mE|Y^{n,\eps}_{t}-\bar Y^n_{t}|^2+C/n^{\beta-2},
\end{align*}
where the last inequality follows by Chebyshev's inequality and condition {\bf (H$^M$)}.
Letting $\eps\to0$ first and then $n\to \infty$, we can get the desired result.
%%% $\tau_n^\eps$ can be controlled by using the assumption {\bf (H$^M$)}, that is how we use the truncation technique.
%%% the right hand of the last inequality does not depends on $\tau_n^\eps$ any more.
\end{proof}

\section{Weak convergence with order $(\alpha/2)\wedge1$}

Now we study the weak convergence of the multi-scale system (\ref{sde0}) to the effective system (\ref{sde1}) in the fully coupled case, i.e., the diffusion coefficient $G(t,x,y)$ in the slow part also depends on the fast term. We first prove the following regularity result for the averaged coefficients.

\bl\label{hfff}
Assume that $a, b\in C_b^{\delta,\alpha}$ and $F,G\in C_b^{\alpha/2,\delta,\alpha}$ with $0<\delta,\alpha\leq 2$. Let $\bar F$ and $\bar G$ be defined by (\ref{bf}). Then we have $\bar F, \bar H\in C_b^{\alpha/2,\alpha}$.
\el
\begin{proof}
We only sketch the proof of the regularity for $\bar F$. Note that when  $0<\alpha\leq 1$, the conclusion has been proven in Lemma \ref{hff}. Let us focus on the case $1<\alpha\leq 2$.	We write for $y_1,y_2\in\mR^{d_2}$
	\begin{align*}
	\nabla_y\bar F(t,y_1)-\nabla_y\bar F(t,y_2)&=\int_{\mR^{d_1}}\![\nabla_yF(t,x,y_1)-\nabla_yF(t,x,y_2)]\mu^{y_1}(\dif x)\\
	&\quad+\int_{\mR^{d_1}}F(t,x,y_2)\nabla_y\big[\mu^{y_1}(\dif x)-\mu^{y_2}(\dif x)\big]=:\tilde\sK_1+\tilde\sK_2.
	\end{align*}
	It is easy to see that there exists a constant $C_1>0$ such that
	$$
	\tilde\sK_1\leq C_1\big(|y_1-y_2|^{\alpha-1}\wedge1\big).
	$$
	For the second term, by the same argument as before we write
	\begin{align*}
	\sK_2&=\int_{\mR^d}F(t,x',y)\big[\nabla_yp_\infty(x';y_1)-\nabla_yp_\infty(x';y_2)\big]\dif x'\\
	&=\lim_{t\to\infty}\int_0^t\!\!\int_{\mR^d}p_{t-s}(x,z;y)\bigg[\frac{\p\sL_0(z,y_1)}{\p y}-\frac{\p\sL_0(z,y_2)}{\p y}\bigg]\\
	&\qquad\qquad\qquad\qquad\qquad\qquad\times\left(\int_{\mR^d}p_s(z,x';y)F(t,x',y)\dif x'\right)\dif z\dif s\\
	&=\int_0^\infty\!\!\int_{\mR^d}p_\infty(z;y)\bigg[\frac{\p\sL_0(z,y_1)}{\p y}-\frac{\p\sL_0(z,y_2)}{\p y}\bigg]\left(\int_{\mR^d}p_s(z,x';y)F(t,x',y)\dif x'\right)\dif z\dif s.
	\end{align*}
	Then, the desired estimates follow by exactly the same arguments as in the proof of Lemma \ref{hff}. We omit the details.
\end{proof}

Recall that $\bar\sL$ is defined by (\ref{lby}). Given a function $\varphi\in C_b^{2+\alpha}$ and $T>0$, we consider the following Cauchy problem:
\begin{equation}\left\{\begin{array}{l}\label{PDE}
\displaystyle
\p_t\hat u(t,y)-\bar\sL\hat u(t,y)=0,\quad t\in [0, T),\\
\hat u(0,y)=\varphi(y).
\end{array}\right.
\end{equation}
It is known that there exists a unique  solution $\hat u$ to (\ref{PDE}) which is given by
$$
\hat u(t,y)=\mE\varphi(\bar Y_t(y)).
$$
Moreover, we have $\nabla^2_y\hat u\in C_b^{\alpha/2,\alpha}$, see e.g. \cite[Chapter IV, Section 5]{La-So-Ur}. Set
$$
\tilde u(t,y):=\hat u(T-t,y),\quad t\in [0, T].
$$
By It\^o's formula, we deduce that
\begin{align*}
\tilde u(T,Y_T^\eps)=\tilde u(0,y)+\int_0^T\p_s\tilde u(s,Y_s^\eps)+\sL_1(X_s^\eps,Y_s^\eps)\tilde u(s,Y_s^\eps)\dif s+\tilde M_t,
\end{align*}
where $\tilde M_t$ is a martingale given by
$$
\tilde M_t:=\int_0^tG(s,X_s^\eps,Y_s^\eps)\nabla_y\tilde u(s,Y_s^\eps)\dif W_s^2.
$$
Note that $$
\tilde u(T,Y_T^\eps)=\hat u(0,Y_T^\eps)=\varphi(Y_T^\eps), \quad\text{and}\quad\tilde u(0,y)=\hat u(T,y)=\mE[\varphi(Y_T)],
$$
and
\begin{align*}
&\p_s\tilde u(s,Y_s^\eps)+\sL_1(X_s^\eps,Y_s^\eps)\tilde u(s,Y_s^\eps)=\sL_1(X_s^\eps,Y_s^\eps)\tilde u(s,Y_s^\eps)-\bar\sL_y\tilde u(s,Y_s^\eps)\\
&=[H(s,X_s^\eps,Y_s^\eps)-\bar H(s,Y_s^\eps)]\nabla^2_y\tilde u(s,Y_s^\eps)+[F(s,X_s^\eps,Y_s^\eps)-\bar F(s,Y_s^\eps)]\nabla_y\tilde u(s,Y_s^\eps).
\end{align*}
We thus get
\begin{align}
\mE[\varphi(Y_T^\eps)]&-\mE[\varphi(Y_T)]=\mE\bigg(\int_0^T[H(s,X_s^\eps,Y_s^\eps)-\bar H(s,Y_s^\eps)]\nabla^2_y\tilde u(s,Y_s^\eps)\dif s\bigg)\no\\
&\quad+\mE\bigg(\int_0^T[F(s,X_s^\eps,Y_s^\eps)-\bar F(s,Y_s^\eps)]\nabla_y\tilde u(s,Y_s^\eps)\dif s\bigg):=\sU_1+\sU_2.\label{wea}
\end{align}
Define
$$
\tilde H(t,x,y):=[H(t,x,y)-\bar H(t,y)]\nabla^2_y\tilde u(t,y)
$$
and
$$
\tilde F(t,x,y):=[F(t,x,y)-\bar F(t,y)]\nabla_y\tilde u(t,y).
$$
Let $\tilde H_n$, $\tilde F_n$ be the mollifying approximations of $\tilde H$ and $\tilde F$ defined similarly as in (\ref{fn}), respectively.
We prepare the following approximation result, which is similar to Lemma \ref{momo}.

\bl
Assume that $a,b\in C_b^{\delta,\alpha}$ and $F\in C_b^{\alpha/2,\delta,\alpha}$ with $0<\delta\leq 1$, $0<\alpha\leq 2$. Then we have
\begin{align}\label{n111}
\|\tilde F-\tilde F_n\|_\infty+\|\tilde H-\tilde H_n\|_\infty+\|a-a_n\|_\infty+\|b-b_n\|_\infty\leq C_0n^{-\alpha},
\end{align}
and
\begin{align}\label{n222}
\|\tilde F_n\|_{C_b^{1,\delta,\alpha}}+\|\tilde F_n\|_{C_b^{\alpha,\delta,2}}+\|\tilde H_n\|_{C_b^{1,\delta,\alpha}}+\|\tilde H_n\|_{C_b^{\alpha,\delta,2}}+ \|a_n\|_{C_b^{\delta,2}}+\|b_n\|_{C_b^{\delta,2}}\leq C_0n^{2-\alpha},
\end{align}
where $C_0>0$ is a constant independent of $n$.
\el
\begin{proof}
Note that when $0<\alpha\leq 1$, the conclusion has been proved in Lemma \ref{hff}. Below, we shall focus on the case $1<\alpha\leq 2$, and prove the corresponding estimates for $\tilde H$. The other estimates can be proved similarly.	According to Lemma \ref{hfff}, it is easy to check that $\tilde H\in C_b^{\alpha/2,\delta,\alpha}$. As a result, we have
\begin{align*}
|\tilde H(t,x,y)-\tilde H_n(t,x,y)|&\leq \int_{\mR^{d_2+1}}\big|\tilde H(t-s,x,y-z)+\tilde H(t-s,x,y+z)\\
&\quad\quad\qquad\quad-2\tilde H(t,x,y)\big|\cdot\rho_2^{n}(z)\rho_1^n(s)\dif z\dif s\\
&\leq C_1\int_{\mR^{d_2+1}}\big(s^{\alpha/2}+|z|^{\alpha}\big)\cdot\rho_2^{n}(z)\rho_1^n(s)\dif z\dif s\leq C_1n^{-\alpha},
\end{align*}
and
\begin{align*}
|\nabla_y^2\tilde H_n(t,x,y)|&\leq \int_{\mR^{d_2+1}}\big|\nabla_y\tilde H(t-s,x,y-z)-\nabla_y\tilde H(t-s,x,y)\big|\cdot|\nabla_z\rho_2^{n}(z)|\rho_1^n(s)\dif z\dif s\\
&\leq C_2n\int_{\mR^{d_2+1}}|z|^{\alpha-1}\cdot\rho_2^{n}(z)\rho_1^n(s)\dif z\dif s\leq C_2n^{2-\alpha}.
\end{align*}
So, the proof is finished.
\end{proof}

We are now in the position to give:

\begin{proof}[Proof of Theorem \ref{main2}]
We begin from (\ref{wea}) and proceed to control the first term. We write
\begin{align*}
\sU_1&\leq \mE\left|\int_0^T\big[\tilde H(s,X_s^\eps,Y_s^\eps)-\tilde H_n(s,X_s^\eps,Y_s^\eps)\big]\dif s\right|\\
&\quad+\mE\left(\int^{T}_0\tilde H_n(s,X^{\eps}_s, Y^{\eps}_s)\dif s\right)=:\sU_{11}+\sU_{12}.
\end{align*}
Using (\ref{n111}), we can control the first term easily by
$$
\sU_{11}\leq C_1n^{-\alpha}.
$$
To control the second term, let $\tilde\Psi_n$ be the solution to the following Poisson equation in $\mR^{d_1}$:
$$
\sL_0^n(x,y)\tilde\Psi_n(t,x,y)=\tilde H_n(t,x,y),
$$
where $\sL_0^n$ is defined by (\ref{l0n}) and $(t,y)\in\mR_+\times\mR^{d_2}$ are viewed as parameters.
Note that $\tilde H_n$ satisfies the centering condition (\ref{cen}). Thus, according to Theorem \ref{popde},
we can use the It\^o's formula to get that
\begin{align*}
\mE\left(\int^{T}_0\tilde H_n(s,X^{\eps}_s, Y^{\eps}_s)\dif s\right)
&=\eps\tilde\Psi_n(T,X_{T}^\eps,Y_{T}^\eps)-\eps\tilde\Psi_n(0,x,y)\\
&\quad+\bigg(\int_0^{T}\big[b_n(X_s^\eps,Y_s^\eps)-b(X_s^\eps,Y_s^\eps)\big]\nabla_x\tilde\Psi_n(s,X_s^\eps,Y_s^\eps)\dif s\\
&\quad+\int_0^{T}\big[a_n(X_s^\eps,Y_s^\eps)-a(X_s^\eps,Y_s^\eps)\big]\nabla^2_x\tilde\Psi_n(s,X_s^\eps,Y_s^\eps)\dif s\bigg)\\
&\quad-\eps\int_0^{T}\big(\p_s+\sL_1\big)\tilde\Psi_n(s,X_s^\eps,Y_s^\eps)\dif s\\
&=:\tilde\sQ_1(T,\eps)+\tilde\sQ_2(T,\eps)+\tilde\sQ_3(T,\eps).
\end{align*}
Using (\ref{n111}), (\ref{n222}) and exactly the same arguments as before, we  get
$$
\tilde\sQ_1(T,\eps)+\tilde\sQ_2(T,\eps)\leq C(\eps+n^{-\alpha}).
$$
and
$$
\tilde\sQ_3(T,\eps)\leq C\eps n^{2-\alpha}.
$$
As a result, we have
$$
\sU_1\leq C(\eps+n^{-\alpha}+\eps n^{2-\alpha}).
$$
Using exactly the same arguments as above, we can also get
$$
\sU_2\leq C(\eps+n^{-\alpha}+\eps n^{2-\alpha}).
$$
Hence, taking $n=\eps^{-1/2}$, we arrive at
$$
\big|\mE[\varphi(Y_T^\eps)]-\mE[\varphi(\bar{Y}_T)]\big|\leq C_T\big(\eps^{\alpha/2}+\eps\big)\leq C_T\eps^{(\alpha/2)\wedge1}.
$$
The proof is finished.
\end{proof}

Finally, we give:

\begin{proof}[Proof of Theorem \ref{main3}]
It is well-known that the solution $u^\eps$ to equation (\ref{pde22}) has the following probabilistic representation (see \cite{Kr1}):
$$
u^\eps(t,x,y)=\mE\left(\int_0^T\psi(Y_s^\eps)\dif s+\varphi(Y_{T-t}^\eps)\right).
$$
Since $\varphi$ is continuous, we can always find a sequence of functions $\varphi_n\in C_b^{3}$ such that
$\|\varphi_n-\varphi\|_\infty\to0$ as $n\to\infty$. As a result, we  deduce by Theorem \ref{main2} that
\begin{align*}
\mE\varphi(Y_{T-t}^\eps)-\mE\varphi(\bar Y_{T-t})\leq \big[\mE\varphi_n(Y_{T-t}^\eps)-\mE\varphi_n(\bar Y_{T-t})\big]+2\|\varphi_n-\varphi\|_\infty.
\end{align*}
Taking $\eps\to0$ first and then $n\to\infty$, we get
$$
\lim_{\eps\to0}\big|\mE\varphi(Y_{T-t}^\eps)-\mE\varphi(\bar Y_{T-t})\big|=0.
$$	
On the other hand, since $\psi$ is bounded, we can always find a sequence of functions $\psi_n\in C_b^{3}$ such that for every $p\geq 1$,
$\|\psi_n-\psi\|_{L^p_{loc}}\to0$ as $n\to\infty$. Then, for every $R>0$, we write
\begin{align*}
&\mE\left(\int_0^T\big[\psi(Y_s^\eps)-\psi(\bar Y_s)\big]\dif s\right)=\mE\left(\int_0^T\big[\psi_n(Y_s^\eps)-\psi_n(\bar Y_s)\big]\dif s\right)\\
&+\mE\left(\int_0^T\big[\psi(Y_s^\eps)-\psi_n(Y_s^\eps)\big]1_{\{|Y_s^\eps|\leq R\}}\dif s\right)+\mE\left(\int_0^T\big[\psi(\bar Y_s)-\psi_n(\bar Y_s)\big]1_{\{|\bar Y_s|\leq R\}}\dif s\right)\\
&+\mE\left(\int_0^T\big[\psi(\bar Y_s)-\psi_n(\bar Y_s)\big]1_{\{|\bar Y_s|> R\}}\dif s\right)+\mE\left(\int_0^T\big[\psi(\bar Y_s)-\psi_n(\bar Y_s)\big]1_{\{|\bar Y_s|> R\}}\dif s\right).
\end{align*}
Due to Theorem \ref{main2}, the first term goes to $0$ as $\eps\to0$.
By Krylov's estimate (see \cite{Kr1}) we have that for some $p>d_1+d_2$,
\begin{align*}
\mE\left(\int_0^T\big[\psi(Y_s^\eps)-\psi_n(Y_s^\eps)\big]1_{\{|Y_s^\eps|\leq R\}}\dif s\right)\leq C\|\psi-\psi_n\|_{L^p_{loc}},
\end{align*}
which goes to 0 as $n\to\infty$. Finally, the last part goes to 0 as $R\to\infty$ by Chebyshev's inequality. This finishes the proof.
\end{proof}

\end{document}